\documentclass[]{my_class}

\SHORTTITLE{Stability of Noisy Tilings}

\TITLE{On the Besicovitch-Stability\\of Noisy Random Tilings\thanks{
This work was supported by the
ANR ``Difference'' project (ANR-20-CE40-0002) and the CIMI LabEx
``Computability of asymptotic properties of dynamical systems'' project
(ANR-11-LABX-0040).

The second author wishes to thank Siamak Taati
for many fruitful references on the stability of tilings.}} % Optional

\AUTHORS{Gayral~Léo\footnote{
University Toulouse III - Paul Sabatier, France.
\BEMAIL{lgayral@math.univ-toulouse.fr}
\url{https://perso.ens-lyon.fr/leo.gayral/en/}}
\and
Sablik~Mathieu\footnote{
University Toulouse III - Paul Sabatier, France
\BEMAIL{msablik@math.univ-toulouse.fr}
\url{https://www.math.univ-toulouse.fr/\~msablik/index.html}}
}


\AMSSUBJ{37B51;37A50} % Separate items with ;
\AMSSUBJSECONDARY{60K35;82B43} % Optional, separate items with ;

\SUBMITTED{January 1, 2021} % Edit.
\ACCEPTED{January 2, 3021} % Edit.

\VOLUME{0}
\YEAR{2012}
\PAPERNUM{0}
\DOI{vVOL-PID}

\newcommand{\N}{\mathds{N}}
\newcommand{\Z}{\mathds{Z}}

\newcommand{\A}{\mathcal{A}}
\newcommand{\B}{\mathcal{B}}
\newcommand{\F}{\mathcal{F}}
\newcommand{\G}{\mathcal{G}}
\newcommand{\M}{\mathcal{M}}

\renewcommand{\epsilon}{\varepsilon}
\renewcommand{\phi}{\varphi}

\usepackage{stmaryrd}

\usepackage{tikz}
\newlength{\mytextsize}
\newcommand{\rbtile}{\resizebox{.7\mytextsize}{!}{\begin{tikzpicture}
    \draw [line width = .5\mytextsize] (-1,-1.09) -- (-1,1.09);
    \draw [line width = .5\mytextsize] (-1,-1) -- (1,-1);
    \draw [line width = .5\mytextsize] (1,-1.09) -- (1,1.09);
    \draw [line width = .5\mytextsize] (-1,1) -- (1,1);
    \draw [line width = .5\mytextsize] (0,0.09) -- (0,-1);
    \draw [line width = .5\mytextsize] (0,0) -- (1,0);
\end{tikzpicture}}}
\newcommand{\lbtile}{\resizebox{.7\mytextsize}{!}{\begin{tikzpicture}
    \draw [line width = .5\mytextsize] (-1,-1.09) -- (-1,1.09);
    \draw [line width = .5\mytextsize] (-1,-1) -- (1,-1);
    \draw [line width = .5\mytextsize] (1,-1.09) -- (1,1.09);
    \draw [line width = .5\mytextsize] (-1,1) -- (1,1);
    \draw [line width = .5\mytextsize] (0,0.09) -- (0,-1);
    \draw [line width = .5\mytextsize] (0,0) -- (-1,0);
\end{tikzpicture}}}
\newcommand{\rttile}{\resizebox{.7\mytextsize}{!}{\begin{tikzpicture}
    \draw [line width = .5\mytextsize] (-1,-1.09) -- (-1,1.09);
    \draw [line width = .5\mytextsize] (-1,-1) -- (1,-1);
    \draw [line width = .5\mytextsize] (1,-1.09) -- (1,1.09);
    \draw [line width = .5\mytextsize] (-1,1) -- (1,1);
    \draw [line width = .5\mytextsize] (0,0.09) -- (0,1);
    \draw [line width = .5\mytextsize] (0,0) -- (1,0);
\end{tikzpicture}}}
\newcommand{\lttile}{\resizebox{.7\mytextsize}{!}{\begin{tikzpicture}
    \draw [line width = .5\mytextsize] (-1,-1.09) -- (-1,1.09);
    \draw [line width = .5\mytextsize] (-1,-1) -- (1,-1);
    \draw [line width = .5\mytextsize] (1,-1.09) -- (1,1.09);
    \draw [line width = .5\mytextsize] (-1,1) -- (1,1);
    \draw [line width = .5\mytextsize] (0,0.09) -- (0,1);
    \draw [line width = .5\mytextsize] (0,0) -- (-1,0);
\end{tikzpicture}}}

\usepackage{float}

\ABSTRACT{
In this paper, we introduce a noisy framework for SFTs,
allowing some amount of forbidden patterns to appear.
Using the Besicovitch distance,
which permits a global comparison of configurations,
we then study the closeness of noisy measures to non-noisy ones
as the amount of noise goes to $0$.
Our first main result is the full classification of the (in)stability
in the one-dimensional case.
Our second main result is a stability property under Bernoulli noise
for higher-dimensional periodic SFTs,
which we finally extend to an aperiodic example
through a variant of the Robinson tiling.
}

\begin{document}

\section{Introduction}

A subshift of finite type (usually called a \emph{SFT})
is the dynamical version of a tiling defined by local rules.
Given a finite set of forbidden patterns on the finite alphabet $\A$,
the corresponding SFT is the set of $\A$-colourings of $\Z^d$
where no forbidden pattern appears.
In the last decades, there has been numerous studies on how
local constraints affect the global structure of tilings in non-trivial ways.
The most important property studied is certainly aperiodicity:
there exists some local rules which impose that any tiling (or configuration)
of the SFT is not periodic~\cite{Ber66,Rob71,Kari96}.
Such SFTs are said to be \emph{aperiodic}.

Aperiodic SFTs received a strong wave of interest when non-periodic crystals
(now called \emph{quasicrystals}) were discovered by
the chemists Dan Shechtman \emph{et al.}~\cite{SheBleGraCa84}.
The connection with tilings was indeed quickly made,
with tiles representing atom clusters,
and forbidden patterns modelling constraints
on the way these atoms can fit together,
\emph{e.g.}\ finite range energetic interactions between them~\cite{LeviStein84}.
In computer science,
tilings have been used as static geometrical models of computation,
ever since Berger proved the undecidability of the so-called \emph{domino problem}
by implementing Turing machines in aperiodic tilings~\cite{Ber66}.
More recently, simulations of Turing machines
have been implemented in several different ways
in order to construct \emph{complex} tilings~\cite{Hoch09,AuSab10,DuRoShe12}.
Yet, in each construction, the aperiodic structure of the tiling is
the key element to embed computations.

A natural question is whether aperiodic structures in SFTs
survive in the presence of some amount of noise,
considering how real-life quasicrystals can have some defects.
In other words, we want to know if a configuration with few mistakes
is structurally close to a generic configuration of the SFT
up to a small amount of changes (ideally proportional to the number of mistakes).
A first step in this direction is the construction
of a three-dimensional model with infinite-range interactions which,
at low positive temperature,
enforces the Thue-Morse sequence along one direction~\cite{AerMiZa98}.
Empirical evidence, obtained through Monte Carlo simulations,
also suggests that even in two dimensions,
using Ammann's aperiodic tileset (famous for using only 16 tiles),
valid tilings remain stable at sufficiently low positive temperatures~\cite{AriRad11}.
In a more formal way, there is a four-dimensional model
with finite-range interactions which, at low but positive temperatures,
admits Gibbs measures that are perturbations of Ammann's aperiodic tiling
along two directions~\cite{Taati}.

The question is also natural in the context of tilings as static models of computation.
Indeed, for a given model of computation, 
we want to know whether the model is robust to errors,
as in the case of cellular automata~\cite{Ga01} or Turing machines~\cite{AsaCol05}.
We consider here $\epsilon$-Bernoulli noises, such that all the cells
may violate local rules with probability $\epsilon$,
independently from each other.
In the case of general tilings,
fixed-point methods allow us to construct an aperiodic SFT
that is robust in the presence of Bernoulli noises~\cite{DuRoShe12}.
In that case, robustness means that for any proportion $\alpha>0$,
for small-enough values of $\epsilon>0$,
any typical noisy configuration is close to some valid configuration of the SFT,
up to a subset of cells of density at most $\alpha$ in $\Z^2$.
The key of the proof is to obtain a tiling where we can repair
\emph{islands} of errors, delimited areas containing forbidden patterns
with a large neighbourhood without any other violation of local rules.
Such matters will be discussed again in Section~\ref{sec:RobustDurRoShe}.

In this paper, we provide a quantitative formalism of error robustness
for a given choice of local rules.
In other words, we allow a small proportion of cells of $\Z^d$ to violate local rules,
and we want to quantify how close a generic noisy configuration is
to a generic non-noisy configuration.
To do so, it is easier to use a distance on the associated measure spaces.
Formally, we will say that a SFT is $f$-stable for a given distance
on the probability measure space if any shift-invariant measure
with a proportion (at most) $\epsilon$ of errors
is at distance at most $f(\epsilon)$ of the set of
non-noisy invariant measures of the SFT.
This formalism will be precisely introduced in Section~\ref{sec:Framework},
up to the notion of stability that begins Section~\ref{sec:Besicovitch}.

Closeness in the weak-* topology is not enough,
because it only characterises a high probability of agreement
on a large \emph{finite} box around the origin,
but does not say anything about the symbols
that are arbitrarily far from the origin.
Thus, when the proportion of errors goes to zero,
the shift-invariant noisy measures must
necessarily converge to non-noisy measures of the SFT,
which we demonstrate in Subsection~\ref{subsec:weakStability}.
Yet, we can exhibit generic noisy configurations,
with an arbitrarily small proportion of errors $\epsilon$,
which we cannot superpose with any configuration of the SFT
on a high-density subset of $\Z^d$.

In order to compare configurations of $\A^{\Z^d}$ in a global way,
the Hamming-Besicovitch pseudo-distance~\cite{BlaForKur97}
-- the density of the set of differences between two configurations in $\Z^d$ --
is a natural approach, as it gives the same importance to all cells.
It is possible to transpose this pseudo-distance into a genuine distance
on the set of measures,
in a similar fashion to the Kantorovich metric~\cite{Ver06}.
This process is detailed in Subsection~\ref{subsec:Besicovitch}.

In this paper, after introducing the general framework of stability,
we prove in Subsection~\ref{subsec:Conjugacy}
that this notion is an invariant of conjugacy.
Thus the stability of a tiling does not depend on the local rules used to define it.
In Section~\ref{sec:1D} we characterise which one-dimensional SFTs
are stable (Theorem~\ref{thm:Characterisation1D}).
It is well known that a one-dimensional SFT is represented
as the set of bi-infinite paths in a labelled transition graph~\cite{LiMa21},
also called an \emph{automaton}.
If this automaton has an aperiodic structure,
then it is possible to correct linearly the mistakes.
On the contrary,
if it has a periodic structure, it is impossible to repair a mistake
at the interface between two phases misaligned within the automaton.
These one-dimensional (un)stable examples can then be extended
to (un)stable SFTs in any dimension with Corollary~\ref{cor:ExtensionDimension}.

Unlike the one-dimensional case,
we prove in Section~\ref{sec:2DPeriodic} that
bi-dimensional strongly periodic SFTs are
linearly stable (Theorem~\ref{thm:PeriodicStability}).
The main idea here is that, if local rules are respected on some region of $\Z^2$,
then this region is the restriction of a periodic configuration,
except maybe on the boundary of the region.
A percolation argument then allows us to prove uniqueness of
such an \emph{infinite} region, what's more with a linear control on its density.
This linear $O(\epsilon)$-stability result is to put into perspective with 
the $O\left(1/\sqrt{\ln(1/\epsilon)}\right)$-stability
obtained in Section~\ref{sec:RobustDurRoShe},
using the strategy described in
\emph{Fixed-point tile sets and their applications}~\cite{DuRoShe12},
which holds in particular for periodic tilings~\cite{BaDuJean10}.

In Section~\ref{sec:Robinson} we consider the famous
Robinson aperiodic tiling~\cite{Rob71} and,
up to some modifications, we show that it is
$O\left(\sqrt[3]{\epsilon}\right)$-stable (Theorem~\ref{thm:RobinsonStability}).
The key idea here is that Robinson configurations are almost periodic,
up to a low-density grid of cells,
so that we may adapt the periodic percolation argument
from Section~\ref{sec:2DPeriodic}.
This result is interesting for two reasons.
First, the Robinson tiling is not robust
in the sense of Durand, Romashchenko and Shen~\cite{DuRoShe12},
so it provides a new, perhaps simpler example of stable aperiodic tiling.
Second, the speed obtained here, though not linear,
is still polynomial,
so much faster than the one one from Section~\ref{sec:RobustDurRoShe}.
The question of whether we can achieve linear stability
for some aperiodic SFT remains open.

\section{Noisy Framework and Weak-* Stability} \label{sec:Framework}

\subsection{Noisy Framework}

\begin{definition}[Configuration Space]
Consider the network $\Z^d$ with $d\in\N^*$ a positive integer,
and a finite alphabet $\A$.
The full-shift configuration space is $\Omega_\A=\A^{\Z^d}$.

We endow this space with the discrete product topology.
In this framework, the clopen cylinders $[w]=\left\{\omega\in\Omega_\A,
\omega|_I = w\right\}$ form a countable base of the topology,
where $w\in \A^I$ is a \emph{finite pattern} over a
\emph{window} $I\varsubsetneq \Z^d$, a \emph{finite} ranges of cells.
Consequently, we use the induced Borel algebra on this space.

For any vector $k\in\Z^d$, let us define the shift
$\sigma_k:\Omega_\A\to \Omega_\A$
such that $\sigma_k(\omega)_l = \omega_{k+l}$.
Likewise, we can apply $\sigma_k$ to any finite pattern or even a range of cells.

Let us denote $\left(e_i\right)_{1\leq i \leq d}$ the canonical basis of $\Z^d$.
Then $\left(\Omega_\A,\sigma_{e_1},\dots,\sigma_{e_d}\right)$ forms
a commutative dynamical system.
\end{definition}

\begin{definition}[Subshift of Finite Type]
A \emph{subshift} of $\Omega_\A$ is a $\sigma$-invariant subset,
\emph{i.e.}\ a subset stable under the action of any shift $\sigma_k$ with $k\in \Z^d$.

Let $\F$ be a finite set of \emph{forbidden} (finite) patterns $w\in \A^{I(w)}$.
A SFT is the subshift $\Omega_\F$ induced by such a set $\F$ as follows:
\[
\Omega_\F:=\left\{ \omega\in \Omega_\A ,\forall w \in\F,
\forall k\in\Z^d,\sigma_k(\omega)|_{I(w)}\neq w \right\} .
\]
In other words, the configurations of the SFT are the ones
where no forbidden pattern occurs anywhere.
Then $\left(\Omega_\F,\sigma_{e_1},\dots,\sigma_{e_d}\right)$
is a commutative dynamical system.
\end{definition}

For a measure $\mu$ on $\Omega$
and a measurable mapping $\theta:\Omega\to\Omega'$,
we denote $\theta^*(\mu)$ the pushforward measure on $\Omega'$
such that $\left[\theta^*(\mu)\right](B) = \mu\left(\theta^{-1}(B)\right)$
for any measurable set $B$.

\begin{definition}[Invariant Probability Measures]
A measure $\mu$ is said to be $\sigma$-invariant if, for any $k\in\Z^d$,
$\sigma_k^*(\mu)=\mu\circ \sigma_{-k}$ is equal to $\mu$.
For a given SFT induced by $\F$,
we denote $\mathcal{M}_\F$ the set of $\sigma$-invariant probability measures
on the space $\Omega_\F$ --
and $\M_\A$ for all the $\sigma$-invariant probability measures
on the full-shift $\Omega_\A$.
\end{definition}

By compactness, this set is non-empty as long as $\Omega_\F\neq \emptyset$.
We will always work under that assumption onwards.
Let us now introduce our noisy \emph{clair-obscur} framework:

\begin{definition}[Noisy SFT]
Consider the alphabet $\widetilde{\A}=\A\times \left\{0,1\right\}$.
For a given cell value $(a,b)\in\widetilde{\A}$,
whenever $b=0$ we say that the cell is \emph{clear},
and when $b=1$ we say that the cell is \emph{obscured}.
We may thus identify $\A$ to the clear subset
$\A\times \left\{0\right\} \varsubsetneq \widetilde{\A}$.
Formally, we will denote $\pi_1:\widetilde{\A}\to \A$
and $\pi_2:\widetilde{\A}\to\left\{0,1\right\}$ the canonical projections.

By extension, we will call patterns $w\in\widetilde{\A}^I$
and configurations $\omega\in\Omega_{\widetilde{\A}}$
\emph{clear} when they are actually defined on the alphabet $\A$,
in opposition to the \emph{obscured} ones,
that contain at least one \emph{obscured} letter in $\widetilde{\A}\backslash \A$.

Using the same identification,
we can define the set of forbidden clear patterns as $\widetilde{\F}
=\left\{ \left(w, 0^{I(w)}\right)\in\widetilde{A}^{I(w)}, w\in\F\right\}$
and the corresponding SFT
on the space $\Omega_{\widetilde{\F}} \subset \Omega_{\widetilde{\A}}$.
\end{definition}

\begin{remark}[Noise vs. Impurities]
With the notion of noise defined above,
comparing the \emph{clear} configurations on $\A$
is a mere matter of projecting $\pi_1:\Omega_{\widetilde{\F}}\to \Omega_\A$,
which results in configurations that may have some amount of forbidden patterns.

Another way to define noise would be to add a blank symbol $\square\notin\A$
not already in the alphabet, without changing $\F$.
The main difference in this case is that there is no natural way to project
$\square$ into $\A$ so that we can compare \emph{clear} configurations.
The symbol $\square$ behaves less like a noise and more like an impurity in itself.

From the point of view of the entropy, this changes things up.
Informally, when the binary noise is maximal,
we can obtain the uniform measure on $\Omega_\A$, for which the entropy is maximal.
In comparison, the only measure that maximizes the amount of impurities
is the Dirac measure $\delta_{\square^{\Z^d}}$ which has a null entropy.
Studying more precisely the behaviour of the entropy in either of these settings,
as a function of the amount of noise, may yield interesting further results.
\end{remark}

\begin{definition}[Locally and Globally Admissible]
A pattern or configuration on $\widetilde{\A}$ will be called
\emph{locally admissible} whenever it contains no forbidden
clear pattern from $\widetilde{\F}$.

We purposefully set aside the term \emph{globally admissible} for
\emph{locally admissible clear configurations}
(that belong to $\Omega_\F$), and all the patterns they contain.
Please note that, in a nondescript SFT, a \emph{clear locally admissible pattern}
is not necessarily a \emph{globally admissible} one.
\end{definition}

\begin{remark}[Reconstruction Function] \label{rmk:LocalGlobalGeneric}

Consider $\phi:\mathcal{P}_F\left(\Z^d\right)\to \N$ the reconstruction function
defined on finite windows $I\subset \Z^d$ by:
\[
\phi(I)=\inf\left\{ k\in\N, w\in \A^{I+B_k} \text{ is locally admissible }
\Rightarrow w|_I \text{ is globally admissible}\right\},
\]
with $B_k=\llbracket -k,k\rrbracket^d$ the ball of radius $k\in\N$
(for the $\lVert .\rVert_\infty$ norm), thus a $(2k+1)$-square.

\emph{A priori}\ $\phi(I)$ could be infinite.
However, as $\A^I$ is finite,
consider $L_I\subset \A^I$ the finite subset of locally admissible patterns
that are \emph{not} globally admissible.
If $v\in L_I$ could be embedded into arbitrarily large admissible patterns,
then by compactness it would be globally admissible.
In other words, there exists a rank $k(v)\in \N$ after which
there is no locally admissible $w\in\A^{I+B_k}$ such that $w|_I = v$.
This holds for all the patterns in $L_I$,
so that $\phi(I)=\max_{v\in L_I}k(v)<\infty$.

As we can embed Turing machines into SFTs,
the function $(I,\F)\mapsto\phi_\F(I)$ is a sort of non-computable busy beaver.
We will see later on specific choices of $\F$ for which $\phi$ is bounded.
This function may seem anecdotal at first,
but it will in fact appear in some way or another in most of the next sections,
and be a fundamental tool in all of our main results.
\end{remark}

Just as before, one can consider noisy measures on $\widetilde{\A}$ with the space
$\M_ {\widetilde{\F}}$.
However, by doing so, we have no control on the weight of obscured cells,
which is why we introduce the following measure spaces.

\begin{definition}[Noisy Probability Measures]
Let $\epsilon\in[0,1]$.
A $\sigma$-invariant probability measure $\nu$ on $\left\{0,1\right\}^{\Z^d}$
($\nu\in\M_{\{0,1\}}$) is called an $\epsilon$-noise
if the probability of a given cell being obscured is at most $\nu([1])\leq \epsilon$.

For a given class of noises $\mathcal{N}\subset \M_{\{0,1\}}$,
we now define the measure space:
\[
\widetilde{\M_\F^{\mathcal{N}}}(\epsilon) =\left\{
\lambda \in \M_ {\widetilde{\F}},
\pi_2^*(\lambda)\in\mathcal{N} \text{ is an $\epsilon$-noise}\right\} .
\]
Likewise, we define the projection $\M_\F^\mathcal{N}(\epsilon) =
\pi_1^*\left(\widetilde{M_\F^\mathcal{N}}(\epsilon) \right)$,
which consists of measures on $\Omega_\A$.
If no class is written, it is implied that $\mathcal{N}=\M_{\{0,1\}}$,
that we allow for \emph{any} noise.
\end{definition}

\begin{definition}[Classes of Dependent Noises]
We define $\B=\left\{\B(\epsilon)^{\otimes \Z^d},\epsilon\in [0,1]\right\}$
the class of independent Bernoulli noises,
where each cell is obscured with probability $\epsilon$
independently of the other cells.

More generally, we consider the class $\mathcal{D}_k$ of $k$-dependent noises,
such that any two windows at distance at least $k$ are independent.
More formally, $\nu\in\mathcal{D}_k$ when,
for any patterns $w\in \A^I$ and $w'\in\A^J$ such that $d_\infty(I,J)\geq k$,
$\nu([w]\cap[w'])=\nu([w])\nu([w'])$.
For any rank $k$, we naturally have $\mathcal{D}_k\subset \mathcal{D}_{k+1}$.
In particular, $\mathcal{D}_1=\B$.
\end{definition}

A direct consequence of this definition is that, on any class $\mathcal{N}$,
for $\epsilon < \delta$, we have the increasing inclusion
$\widetilde{\M_\F^\mathcal{N}}(\epsilon) \subset 
\widetilde{\M_\F^\mathcal{N}}(\delta)$, 
which naturally still holds for $\M_\F^\mathcal{N}$ after projection.
Let us notice that $\M_\F(0)=\M_\F$ is non-empty,
and that $\M_\F(1)=\M_\A$ is the set of shit-invariant measures on $\Omega_\A$.

We are now interested in the \emph{stability} of noisy measures,
\emph{i.e.}\ in the fact that $\M_\F(\epsilon)$ gets close to $\M_\F$
in some sense -- for some topology -- as $\epsilon$ goes to $0$.

\subsection[Weak-* Stability]{Weak-$*$ Stability} \label{subsec:weakStability}

A natural topology on measures to consider first is the weak-$\ast$ topology,
but we will see here that as $\epsilon\to 0$,
any adherence point of a sequence of noisy measures is in $\M_\F$.

\begin{definition}[Weak-$\ast$ Topology]
We can define the weak-$\ast$ topology
on the space of probability measures on $\A^{\Z^d}$
as the smaller topology such that,
for any finite pattern $w$,
the evaluation $\mu\mapsto \mu([w])$ is continuous.
\end{definition}
Note that this topological space is Hausdorff-compact,
and that the subset $\M_\A$ of $\sigma$-invariant measures is a closed subset.

\begin{lemma}[]
Consider $\mu\in \M_\A$ a $\sigma$-invariant measure on $\A$.
If, for any pattern $w\in\F$, we have $\mu([w])=0$,
then $\mu\in \M_\F$.

\begin{proof}
To show that $\mu\in \M_\F$,
we need to show that the measure is supported by $\Omega_\F$.
The complementary of $\Omega_\F$ is the set
$\bigcup_{k\in \Z^d} \bigcup_{w\in\F} \sigma_k([w])$.
By $\sigma$-invariance, for any $w\in\F$ and $k\in\Z^d$,
we have $\mu\left(\sigma_k([w])\right)=\mu([w])=0$,
thus $\mu\left(\Omega_\F^c\right)=0$,
so that $\mu$ is indeed supported by $\Omega_\F$.
\end{proof}
\end{lemma}

\begin{proposition} \label{WeakStability}
Let $\mu_n\in \M_\F\left(\epsilon_n\right)$ be a sequence of noisy measures,
with $\epsilon_n\underset{n\to\infty}{\longrightarrow} 0$.
Then any adherence value of the sequence is in $\M_\F$.

\begin{proof}
Consider a weakly converging subsequence $\mu_{\theta(n)}\to^* \mu$,
with $\theta:\N\to\N$ an increasing extraction.
Naturally, the limit $\mu$ is also $\sigma$-invariant.

Notice that, as $\epsilon$-noises are defined by
forcing the measure of the noise cylinder $[1]$
to belong to the closed set $[0,\epsilon]$,
the set $\widetilde{\M_\F}(\epsilon)$ is naturally weakly closed.
Hence, by monotonous inclusion,
for any $\epsilon>0$ we have $\mu\in \M_\F(\epsilon)$.

Consider $\lambda_\epsilon\in \widetilde{\M_\F}(\epsilon)$ that projects to $\mu$,
and a forbidden pattern $w\in \F$.
If $[w]$ occurs for $\mu$,
then at least one of the cells of the window $I(w)$
must be obscured for $\lambda_\epsilon$,
so by union bound, $\mu([w])\leq \lambda_\epsilon([1])\times \left| I(w) \right|\leq
|I(w)|\epsilon$.
As $\epsilon$ goes to $0$, we conclude that $\mu([w])=0$.
Using the previous lemma, $\mu\in \M_\F$.
\end{proof}
\end{proposition}

As all SFTs are weakly stable, this property yields no interesting classification.

In the following section,
we will introduce a general notion of metric stability and convergence speed --
as there is no \emph{canonical} metric associated to the weak topology,
we did not try to quantify the speed of convergence in this case.

The main issue with the weak-$*$ topology
is that it looks at things on a local scale,
on finite patterns, without really forcing any kind of behaviour
on $\Z^d$ as a whole.
To better discriminate between SFTs,
we will introduce the Besicovitch distance $d_B$ on measures,
that looks at configurations globally and quantifies the frequency of differences.

Thereafter,
we will prove that the stability for $d_B$ is conjugacy-invariant,
in order to illustrate how this distance is manipulated.

\section{Stability under Besicovitch Topology and Conjugacy Invariance} \label{sec:Besicovitch}

Before diving into the Besicovitch world, let us briefly introduce our
general notion of stability.

\begin{definition}[Stability]
Consider here a distance $d$ on $\M_\A$, a noise class $\mathcal{N}\subset\M_{\{0,1\}}$,
and a non-decreasing function $f:[0,1]\to \mathds{R}^+$,
right-continuous in $0$ with $f(0)=0$.

The SFT induced by $\F$ is said to be $f$-stable
for the distance $d$ on the class $\mathcal{N}$ if:
\[
\forall \epsilon\in [0,1],\sup\limits_{\mu\in \M_\F^\mathcal{N}(\epsilon)}
d\left(\mu, \M_\F\right) \leq f(\epsilon) .
\]
The SFT is stable if it is $f$-stable for some function $f$.
We say that $\Omega_\F$ is linearly stable (resp. polynomialy stable)
if it is $f$-stable with $f(\epsilon)=O(\epsilon)$
(resp. $f(\epsilon)=O\left(\epsilon^\alpha\right)$ for some $0<\alpha\leq 1$).
\end{definition}

\begin{proposition}
Given a metric $d_*$ which induces the weak-$*$ topology,
any SFT is stable for the distance $d_*$.

\begin{proof}
This result is a consequence of Proposition~\ref{WeakStability}.
Indeed, by contraposition, assume some SFT $\Omega_\F$
is not stable for the distance $d_*$.

Then, there exists a sequence $\epsilon_n\underset{n\to\infty}{\longrightarrow}0$
and measures $\mu_n\in\M_\F\left(\epsilon_n\right)$ such that
$\inf_{n\in\N} d_*\left(\mu_n,\M_\F\right) = d >0$.
By compactness, this sequence admits a weak-$*$ adherence value $\mu$,
and $d_*\left(\mu,\M_\F\right)\geq d$.
In particular, $\mu\notin\M_\F$, which contradicts Proposition~\ref{WeakStability},
hence stability.
\end{proof}
\end{proposition}

Note that this result gives no quantitative bound on the speed of convergence.
In order to obtain such bounds,
we may perhaps make a clever use of the reconstruction function,
but as stated earlier, the lack of canonical choice for $d_*$ discouraged this study.

In this section, we will first introduce the Besicovitch distance $d_B$,
and then prove that stability for $d_B$ is conjugacy-invariant
on the class of all noises $\M_{\{0,1\}}$.
At last, we mention the notion of domination,
which will allow us to extend the conjugacy-invariant stability to $\B$.

\subsection{Besicovitch Topology} \label{subsec:Besicovitch}

In order to compare \emph{measures},
we need first to be able to compare \emph{configurations}.

\begin{definition}[Hamming-Besicovitch distance]
On a finite window $I\subset \Z^d$,
we define the Hamming distance between two finite patterns $x,y\in \A^I$ as:
\[
d_I(x,y)=\frac{1}{|I|} \left| \left\{ k\in I, x_k\neq y_k \right\} \right| .
\]

For a given increasing sequence $\left( I_n\right)_{n\in\N}$
such that $\bigcup\limits_{n\in\N} I_n =\Z^d$,
we can define a pseudometric on $\Omega_\A$, such that for $x,y\in\Omega_\A$:
\[
d_H(x,y)=\limsup\limits_{n\to\infty} d_{I_n}\left( x|_{I_n},y|_{I_n}\right) .
\]
This pseudometric $d_H$ is usually called the Hamming-Besicovitch distance.
Remark that the Hamming distances $d_I$ are clearly measurable for the product topology,
thus so is the limit $d_H$.
\end{definition}

In order to use an extension of Birkhoff's pointwise ergodic theorem,
we need $\left(I_n\right)$ to be a sequence of \emph{boxes} (products of intervals). 
A general statement of this theorem can be found
in \emph{Ergodic Theorems}~\cite[Chapter 6]{KreBru85}.
Let us more specifically use the boxes $B_n=\llbracket -n,n\rrbracket^d$ further on.

\begin{definition}[Besicovitch distance]
A measure $\lambda\in\M_{\A\times \A}$ is said to be a coupling
between the measures $\mu,\nu\in \M_\A$ if $\pi_1^*(\lambda)=\mu$
and $\pi_2^*(\lambda)=\nu$.

For two measure $\mu,\nu\in \M_\A$ we define their Besicovitch distance as:
\[
d_B(\mu,\nu)=\inf\limits_{\lambda\text{ a coupling}}
\int d_H(x,y)\mathrm{d}\lambda(x,y) .
\]
Note that we can always consider the independent coupling $\mu\otimes\nu$,
so the set of couplings is non-empty.
\end{definition}

For the coupling, we can more generally consider
probability measures $\lambda$ on some general space $\Omega$
and two measurable applications $\psi_1:\Omega\to\Omega_\A$
(resp. $\psi_2$) such that $\psi_1^*(\lambda)=\mu$
(resp. $\psi_2^*(\lambda)=\nu$) and consider
$d_H\left(\psi_1(\omega),\psi_2(\omega)\right)\mathrm{d}\lambda(\omega)$
in the integral instead, 
as long as we have $\left(\psi_1,\psi_2\right)^*(\lambda)\in\M_{\A\times\A}$.
We will use this more general version to build couplings
that use additional information,
notably which cells should be obscured in the noisy framework,
or the value of an additional independent random variable.

The Besicovitch distance $d_B$ has been quite used in the recent research literature,
but it was already introduced in earlier works,
sometimes also named $\overline{d}$ as
in \emph{Ergodic Theory via Joinings}~\cite[Chapter 15]{Gla03}.
The main interest of $d_B$, in the context of these works,
is that the measure entropy is continuous for this topology.
Even though $d_B$ has been widely studied,
let us prove here that $d_B$ is indeed a distance,
in order to get acquainted with the notion.

\begin{lemma} \label{lem:BesicovitchDistance}

The function $d_B$ is a distance on $\M_\A$,
and $d_B(\mu,\nu)$
is always \emph{reached} for some coupling between the measures.

\begin{proof}
The function $d_B$ is trivially symmetric,
and $d_B(\mu,\mu)=0$ for any measure $\mu\in\M_\A$.

To prove the triangle inequality,
consider three measures $\mu_1,\mu_2,\mu_3\in\M_\A$.
Consider a coupling $\lambda_{1,2}\in\M_{\A\times\A}$
(resp. $\lambda_{2,3}\in\M_{\A\times\A}$)
between $\mu_1$ and $\mu_2$ (resp. $\mu_2$ and $\mu_3$).
The measures $\lambda_{1,2}$ and $\lambda_{2,3}$ are compatible
in the sense that they share a common projection
$\pi_2^*\left(\lambda_{1,2}\right)= \pi_1^*\left(\lambda_{2,3}\right)=\mu_2$.
Thence, it is known~\cite[Chapter 6]{Gla03} that there exists
a coupling $\lambda_{1,2,3}$ between them,
such that $\left(\pi_1,\pi_2\right)^*\left(\lambda_{1,2,3}\right)=\lambda_{1,2}$
and likewise $\left(\pi_2,\pi_3\right)^*\left(\lambda_{1,2,3}\right)=\lambda_{2,3}$.
In particular, $\lambda_{1,2,3}$ gives us a coupling between $\mu_1$ and $\mu_3$:
\[
\begin{array}{rcl}
d_B\left(\mu_1,\mu_3\right)
&\leq&\int d_H(x,z)\mathrm{d}\lambda_{1,2,3}(x,y,z)\\
&\leq & \int d_H(x,y)+d_H(y,z) \mathrm{d}\lambda_{1,2,3}(x,y,z)\\
&= & \int d_H(x,y)\mathrm{d}\lambda_{1,2}(x,y)
+\int d_H(y,z)\mathrm{d}\lambda_{2,3}(y,z) .
\end{array}
\]
Now, by taking the infimum over all couplings $\lambda_{1,2}$ and $\lambda_{2,3}$
we finally obtain the upper bound $d_B\left(\mu_1,\mu_3\right)\leq
d_B\left(\mu_1,\mu_2\right)+d_B\left(\mu_2,\mu_3\right)$,
the triangle inequality.

Consider now some coupling $\lambda\in\M_{\A\times \A}$
between two measures $\mu,\nu\in\M_\A$.
Using our pointwise ergodic theorem, it follows that
$d_H(x,y)$ is an actual limit $\lambda$-almost-surely,
and that $\int d_H(x,y)\mathrm{d}\lambda(x,y) = \int
\mathds{1}_{\left\{x_0\neq y_0\right\}}\mathrm{d}\lambda(x,y)$.
Hence, it is clear that the mapping
$\lambda\mapsto\int d_H(x,y)\mathrm{d}\lambda(x,y)$
is weakly continuous on $\M_{\A\times\A}$,
so by compactness the distance $d_B$ is reached by some coupling $\lambda$.

Assume now that $d_B(\mu,\nu)=0$ is reached for some coupling $\lambda$.
Then, $\lambda$-almost-surely, we have $x_0=y_0$.
As $\lambda$ is $\sigma$-invariant, it is more generally true for any cell $k\in\Z^d$
that $x_k=y_k$ almost-surely.
By taking the countable intersection of such events, $x=y$ almost-surely,
so $\lambda$ is supported by the diagonal of $\Omega_{\A\times\A}=
\Omega_\A\times \Omega_\A$.
Thence, we have $\mu=\pi_1^*(\lambda)=\pi_2^*(\lambda)=\nu$.
Conversely, distinct measures in $\M_\A$ are distinguishable.
\end{proof}
\end{lemma}

The important part of this lemma is that we need the $\sigma$-invariance
of the measures to conclude that $d_B$ is not a mere
\emph{pseudometric} but a distance,
which shows how the Besicovitch distance is appropriate to our specific setting.

\begin{remark}
If $\mu\in\M_\F(\epsilon)$,
then obscured cells have a frequency $\epsilon$.
For a nondescript set of forbidden patterns $\F$,
it is reasonable to assume that whenever a cell is obscured,
the cell contains a ``wrong'' letter with positive probability,
with respect to any globally admissible configuration $\omega\in \Omega_\F$.
Hence, the best and fastest bounds we can reasonably obtain are linear,
\emph{i.e.}\ $f(\epsilon)=\Omega(\epsilon)$ using the $\Omega$ Landau notation.
\end{remark}

Now that the Besicovitch distance $d_B$ has been properly introduced,
let us prove a that stability for $d_B$ is conjugacy-invariant in some sense.

\subsection{Conjugacy-Invariant Stability} \label{subsec:Conjugacy}

Usually, conjugate SFTs share the same dynamical properties.
Hence, proving the invariance of stability under conjugacy
would imply that stability is indeed a property of a SFT $\Omega_\F$
and not of the specific rules $\F$ used to define it.
In order to study the invariance by conjugacy,
let us first properly define what a conjugacy is.

\begin{definition}[Morphism and Conjugacy]
Consider two sets of forbidden patterns $\F$ on the alphabet $\A_1$
and $\G$ on $\A_2$,
not necessarily on the same alphabet, but on the same grid $\Z^d$.
A morphism from the SFT $\Omega_\F$ to $\Omega_\G$
is a continuous $\sigma$-invariant mapping
$\theta:\Omega_\F\to \Omega_\G$.

Equivalently to $\sigma$-invariance and continuity~\cite{Hed69},
we can define $\theta:\A_1^J\to \A_2$ locally, on a finite window $J\subset \Z^d$,
and extend it on a configuration $x\in\Omega_\F$
so that we have $\theta(x)_k = \theta\left( x|_{J+k}\right)$ for any cell $k\in\Z^d$.

When $\A=\A_1=\A_2$ and $\F=\G=\emptyset$,
$\theta:\Omega_\A\to\Omega_\A$ is also called a \emph{cellular automaton}.

Two SFTs $\Omega_\F$ and $\Omega_\G$ are conjugate
if there is a bijective morphism $\theta:\Omega_\F\to\Omega_\G$,
in which case $\theta^{-1}$ must \emph{also} be a morphism.
\end{definition}

Later on, we will always consider the \emph{local} definition of morphisms
as extensions of local mappings $\theta:\A_1^J\to \A_2$.
An interest of this viewpoint is that any morphism from $\Omega_\F$
to $\Omega_\G$ is actually simply the restriction of
a morphism on the full-shifts $\Omega_{\A_1}$ and $\Omega_{\A_2}$.

\begin{definition}[Thickened Noise] \label{def:thickNoise}

Let $\gamma_n:\{0,1\}^{B_n}\to \{0,1\}$ be the cellular automaton
defined by $\gamma_n(w)=\max_{k\in B_n} w_k$.
We say that $\gamma_n(\omega)$ is $n$-thickened for $\omega\in\Omega_{\{0,1\}}$
in the sense that if the cell $c\in\Z^d$ is obscured in $\omega$,
then its $n$-neighbourhood $c+B_n$ is obscured in $\gamma_n(\omega)$.
\end{definition}

These specific morphisms will allow us to obscure
the forbidden patterns that may appear when using
a morphism or a measurable application on $\Omega_\A$ later on.

\begin{lemma} \label{lem:MorphismStability}

Consider the SFTs $\Omega_\F$ and $\Omega_\G$,
and a morphism from $\Omega_\F$ to $\Omega_\G$,
locally defined as $\theta:\A_1^J\to \A_2$.
Then for any $x,y\in \Omega_{\A_1}$, we have
$d_H\left(\theta(x),\theta(y)\right)\leq D_\theta d_H(x,y)$
with the constant $D_\theta=|J|$.
Consequently, for any measures $\mu,\nu \in \M_{\A_1}$,
we have $d_B\left( \theta^*(\mu),\theta^*(\nu)\right)
\leq D_\theta d_B(\mu,\nu)$.

There exists a radius $r_\theta$ such that the morphism
$\widetilde{\theta}:=\left(\theta,\gamma_{r_\theta}\right):
\Omega_{\widetilde{\A_1}}\to \Omega_{\widetilde{\A_2}}$
satisfies $\widetilde{\theta}\left(\Omega_{\widetilde{\F}}\right)
\subset \Omega_{\widetilde{\G}}$.
Moreover, there is a constant $C_\theta$ such that,
whenever $\gamma_{r_\theta}^*\left(\mathcal{N}\right)\subset \mathcal{N}'$,
for any $\epsilon>0$:
\[
\widetilde{\theta}^*\left(\widetilde{\M_\F^\mathcal{N}}(\epsilon)\right) \subset
\widetilde{\M_\G^{\mathcal{N}'}}\left(C_\theta\times\epsilon\right) .
\]

\begin{proof}
Assume that $\theta(x)_k\neq \theta(y)_k$.
Then $x$ and $y$ \emph{must} differ in at least one cell of the window $J+k$.
Conversely, each cell of $\Z^d$ can appear into \emph{at most} $|J|$ such windows,
so that we naturally obtain the bound
$d_H\left(\theta(x),\theta(y)\right)\leq |J| d_H(x,y)$.
Now, assuming $d_B(\mu,\nu)$ is reached for a coupling
$\lambda\in\M_{\A_1\times\A_1}$,
then $(\theta,\theta)^*(\lambda)\in\M_{\A_2\times\A_2}$
is a coupling between $\theta^*(\mu)$ and $\theta^*(\nu)$,
and we consequently obtain the analogous bound for $d_B$.

Just like $\theta:\A_1^J\to \A_2$ naturally
sends $\Omega_{\A_1}$ onto $\Omega_{\A_2}$,
it sends any finite pattern $v\in \A_1^{J+I}$ onto $\theta(v)\in\A_2^I$.
The ``local'' property that characterises
$\theta\left(\Omega_\F\right)\subset\Omega_\G$
is \emph{not} that it preserves locally admissible patterns,
but that it preserves globally admissible ones.

If a locally admissible pattern $v\in\A_1^{J+I}$ is not globally admissible,
nothing forbids $\theta(v)\in\A_2^I$ from containing forbidden patterns of $\G$.
In such a case, let us extend $v$ into $\omega\in\Omega_{\A_1}$
by filling the empty cells outside of $I+J$ with any letter $a\in\A_1$,
and consider the noise $b=\mathds{1}_{(I+J)^c}$
that obscures all the cells outside of $I+J$.
Then naturally $(\omega,b)\in\Omega_{\widetilde{\F}}$ is locally admissible
but $(\theta(\omega),b)\notin\Omega_{\widetilde{\G}}$ is not.
Thence, we cannot simply extend the morphism
$\theta:\Omega_{\A_1}\to\Omega_{\A_2}$
as $\widetilde{\theta}$ by leaving the second coordinate unchanged.

More precisely, assume that $w=\theta(v)\in\G$ is a forbidden pattern,
with $v\in\A_1^{J+I(w)}$.
Then $v$ must not be a \emph{globally} admissible pattern itself.
As explained in Remark~\ref{rmk:LocalGlobalGeneric},
using the reconstruction function, we have $r(w)=\phi_\F(J+I(w))\in \N$
such that, if we can extend $v$ into a locally admissible pattern
$v\in \A_1^{J+I(w)+B_{r(w)}}$, then $v$ itself must be globally admissible.

Let us define
$r_{\theta}=\max_{w\in\G}r(w)+\max_{c\in J}\lVert c \rVert_\infty$.
Consider $(\omega,b)\in\Omega_{\widetilde{\F}}$.
If $\theta(\omega)$ contains a forbidden pattern $w$ in the window $c+I(w)$,
then it follows that the window $c+J+I(w)$ of $\omega$ is not
globally admissible, so the window $c+I(w)+B_{r_{\theta}}$
of $\omega\in\Omega_{\A_1}$ must not be locally admissible.
As $(\omega,b)$ is locally admissible,
this implies that at least one cell in $c+I(w)+B_{r_\theta}$
must be obscured.
We proved that, if $(\omega,b)\in\Omega_{\widetilde{\F}}$,
then $( \theta(\omega),\gamma_{r_\theta}(b))\in\Omega_{\widetilde{\G}}$,
so $\widetilde{\theta}=(\theta,\gamma_{r_\theta})$ is the morphism we wanted.

Finally, we need to exhibit the constant $C_{\theta}$.
Consider a noisy measure $\lambda\in\widetilde{\M_\F^\mathcal{N}}(\epsilon)$,
with an $\epsilon$-noise $\nu=\pi_2^*(\lambda)\in\mathcal{N}$.
Notice that $\pi_2\circ \widetilde{\theta}=\gamma_{r_\theta}$,
so the noise of $\widetilde{\theta}^*(\lambda)$
is actually $\gamma_{r_\theta}^*(\nu)\in\mathcal{N}'$.
Remark that the clear configuration $0^\infty$ is a fixed point of $\gamma_{r_\theta}$.
Note that, as in the proof of Lemma~\ref{lem:BesicovitchDistance},
using a pointwise ergodic theorem, the amount of noise in $\nu$ is:
\[
\nu([1])=\int \mathds{1}_{\{x_0\neq 0\}} \mathrm{d}\nu(x)
=\int d_H\left(x,0^\infty\right)\mathrm{d}\nu(x) = d_B\left(
\nu,\delta_{0^\infty}\right).
\]
Thus, if we apply the first part of the current lemma
to the current morphism $\gamma_{r_\theta}$, with the set $J=B_{r_{\theta}}$,
we conclude that $d_B\left(\gamma_{r_\theta}^*(\nu),\delta_{0^\infty}\right)
\leq C_{\theta} d_B\left(\nu,\delta_{0^\infty}\right)
\leq C_{\theta}\times \epsilon$ with the constant $C_{\theta}=|J|
=\left( 2 r_{\theta}+1\right)^d$.
At last, $\gamma_{r_\theta}^*(\nu)$ is a $\left(C_{\theta}\epsilon\right)$-noise,
$\widetilde{\theta}^*(\lambda)\in\widetilde{\M_\G^{\mathcal{N}'}}
\left(C_{\theta}\epsilon\right)$.
\end{proof}
\end{lemma}

Assume that the SFT $\Omega_\F$ is $f$-stable,
and that it is sent on $\Omega_\G$ by $\theta$.
Using the lemma, we deduce that the \emph{subset}
$\pi_1^*\left(\widetilde{\theta}^*\left(\widetilde{\M_\F^\mathcal{N}}
(\epsilon)\right)\right)
= \theta^*\left(\M_\F^\mathcal{N}(\epsilon)\right)$
of $\M_\G^{\gamma_{r_{\theta}}(\mathcal{N})}\left(C_{\theta}\epsilon\right)$ is
informally $D_{\theta}\times f$-stable.
However, this still does not give us enough information
to obtain a full-fledged and well-defined stability property for $\G$.
To obtain such a result, we will now assume that $\theta$ is not only a morphism
but a conjugacy between $\Omega_\F$ and $\Omega_\G$.

\begin{theorem}[Conjugacy-Invariant Stability] \label{thm:Conjugacy}

Consider a conjugacy $\theta:\Omega_\F\to\Omega_\G$,
and assume that $\Omega_\F$ is $f$-stable
for $d_B$ on a class $\gamma_{r_{\theta^{-1}}}^*(\mathcal{N})$ of noises.

Then there exists a constant $E$ such that $\Omega_\G$ is $g$-stable on $\mathcal{N}$
with the speed
\[
g:\epsilon\mapsto D_{\theta}f\left(C_{\theta^{-1}}\epsilon\right)+E \epsilon.
\]

\begin{proof}
We will use the result of Lemma~\ref{lem:MorphismStability}
for both $\theta:\Omega_\F\to\Omega_\G$ \emph{and} its inverse
$\theta^{-1}:\Omega_\G\to\Omega_\F$.
Note that, on the larger domain $\Omega_{\A_2}$,
the cellular automaton $\theta\circ \theta^{-1}$ is still well-defined,
but is not necessarily the identity function outside of the domain $\Omega_\G$.
Now, if we consider two measures $\mu,\nu \in \M_{\A_2}$:
\[
\begin{array}{rcl}
d_B(\mu,\nu)&\leq&d_B\left(\mu,\left(\theta\circ\theta^{-1}\right)^*(\mu)\right)\\
& +&d_B\left(\left(\theta\circ \theta^{-1}\right)^*(\mu),
\left(\theta\circ \theta^{-1}\right)^*(\nu)\right)\\
&+& d_B\left(\left(\theta\circ \theta^{-1}\right)^*(\nu),\nu\right) .
\end{array}
\]
The idea behind this back-and-forth is that,
by going from $\Omega_\G$ to $\Omega_\F$,
we reach a stable SFT while still keeping the noise under control,
and then going from $\Omega_\F$ to $\Omega_\G$
allows us to maintain this stability while comparing the new configuration
to the old one.
In particular, if $\nu \in \M_\G$,
then it is supported by the domain $\Omega_\G$,
where $\theta\circ\theta^{-1}$ is the identity function,
so that $d_B\left(\left(\theta\circ \theta^{-1}\right)^*(\nu),\nu\right) =0$.

Consider a measure $\mu\in\M_\G^\mathcal{N}(\epsilon)$,
and $\nu_\F\in\M_\F$ that reaches
$d_B\left(\left(\theta^{-1}\right)^*(\mu),\M_\F\right)$.
If we denote $\nu_\G=\theta^*\left(\nu_\F\right)\in\M_\G$, then:
\[
d_B\left(\mu,\M_\G\right) \leq d_B\left(\mu,\nu_\G\right) 
\leq d_B\left(\mu,\left(\theta\circ\theta^{-1}\right)^*(\mu)\right)
+d_B\left(\left(\theta\circ \theta^{-1}\right)^*(\mu),\nu_\G\right) .
\]

In particular, using Lemma~\ref{lem:MorphismStability} for $\theta^{-1}$,
we know that
$\left(\theta^{-1}\right)^*(\mu)\in\M_\F^{\gamma_{r_{\theta^{-1}}}^*
(\mathcal{N})}\left( C_{\theta^{-1}}\epsilon\right)$.
Thence, using Lemma~\ref{lem:MorphismStability} for $\theta$,
as $\Omega_\F$ is $f$-stable
on $\gamma_{r_{\theta^{-1}}}^*(\mathcal{N})$, we get the bound:
\[
\begin{array}{rcl}
d_B\left(\left(\theta\circ \theta^{-1}\right)^*(\mu),\nu_\G\right)
&\leq & D_{\theta} d_B\left( \left(\theta^{-1}\right)^*(\mu),\nu_\F\right)\\
&=& D_{\theta} d_B\left(\left(\theta^{-1}\right)^*(\mu),\M_\F\right)\\
&\leq& D_{\theta} f\left(C_{\theta^{-1}} \epsilon\right).
\end{array}
\]

To conclude the proof,
we just need to have a linear control on
$d_B\left(\mu,\left(\theta\circ\theta^{-1}\right)^*(\mu)\right)$ as $\epsilon\to 0$.
To do so, we will study $d_H\left(x,\theta\circ\theta^{-1}(x)\right)$
for any $x\in \Omega_{\A_2}$.
More precisely, whenever $(x,b)\in\Omega_{\widetilde{\F}}$, we want a bound
$d_H\left(x,\theta\circ\theta^{-1}(x)\right)\leq E d_H\left(b,0^\infty\right)$.
Assuming such a bound holds,
consider $\lambda\in\widetilde{\M_\G^\mathcal{N}}(\epsilon)$
that projects to $\mu$, which naturally gives a coupling
between $\mu=\pi_1^*(\lambda)$ and $\left(\theta\circ\theta^{-1}\right)^*(\mu)=
\left(\pi_1\circ\widetilde{\theta}\circ\widetilde{\theta^{-1}}\right)^*(\lambda)$.
Then we obtain:
\[
d_B\left(\mu,\left(\theta\circ\theta^{-1}\right)^*(\mu)\right) \leq 
\int\limits_{\Omega_{\widetilde{\G}}}
 d_H\left(x,\theta\circ\theta^{-1}(x) \right) \mathrm{d}\lambda(x,b)
\leq E \int\limits_{\Omega_{\widetilde{\G}}} 
d_H\left(b,0^\infty\right) \mathrm{d}\lambda(x,b) \leq E \epsilon .
\]

The sketch of the proof from now on is pretty much the same
as in Lemma~\ref{lem:MorphismStability}.
Let us suppose that $(x,b)\in\Omega_{\widetilde{\G}}$ and that 
$x_k\neq \theta\circ\theta^{-1}(x)_k$ for some cell $k\in\Z^d$.
Consider the window $J=J_{\theta^{-1}}+J_{\theta}$
such that the value of $\theta\circ\theta^{-1}(x)_k$
only depends on the pattern $x|_{J+k}$.
Let us assume that $0\in J$ without loss of generality.
If $x|_{J+k}$ was globally admissible, then we could extend it
into a globally admissible configuration $y\in\Omega_{\G}$,
such that $\theta\circ\theta^{-1}(x)_k=\theta\circ\theta^{-1}(y)_k
=y_k=x_k$.
This contradicts our hypothesis, so $x_{J+k}$ is not globally admissible.
This means that, using once again the reconstruction function $\phi$
from Remark~\ref{rmk:LocalGlobalGeneric} for the SFT $\Omega_\G$,
$x|_{J+k+B_{\phi(J)}}$ is not locally admissible, so the same windows in $b$
contains at least one obscured cell.
Hence, $d_H\left(x,\theta\circ\theta^{-1}(x)\right)\leq E d_H\left(b,0^\infty\right)$
with the constant $E=\left|J+B_{\phi(J)}\right|$, which concludes the proof.
\end{proof}
\end{theorem}

\begin{corollary}[Conjugacy-Invariance for Stable Noise Classes]
If $\mathcal{N}$ is stable under the action of any $\gamma_n$,
then for any two conjugate SFTs $\Omega_\F$ and $\Omega_\G$,
$\Omega_\F$ is stable (resp. linearly stable, polynomially stable)
on the class $\mathcal{N}$ if and only if $\Omega_\G$ is.
\end{corollary}

In particular, this corollary holds
for the class of all noises $\mathcal{N}=\M_{\{0,1\}}$.
This stability hypothesis is actually quite restrictive.
For example, we naturally have the inclusion $\gamma_n(\B)\subset\mathcal{D}_{2n+1}$
but $\gamma_n(\B)\not\subset\mathcal{D}_{2n}$.
Thence, $\gamma_n(\B)\not\subset \B$,
the previous conjugacy-invariance corollary
does not apply on the class $\mathcal{N}=\B$.

\subsection{Stability and Domination}

We will now introduce the notion of domination between noises,
which will allow us to send $\mathcal{D}_k$ back into $\B$,
in order to obtain a conjugacy-invariant stability result for the class $\B$.

\begin{definition}[Domination]
A Borel set $B\subset\{0,1\}^{\Z^d}$ is said to be increasing if,
for any $b\in B$ and $b'\geq b$ (on each coordinate), we have $b'\in B$.

Consider $\nu_1,\nu_2\in\M_{\{0,1\}}$.
We say that $\nu_2$ dominates $\nu_1$, and we denote $\nu_2\geq \nu_1$,
if $\nu_2(B)\geq \nu_1(B)$ for any increasing Borel set $B$.
Equivalently~\cite[Theorem 2.4]{Lig05},
there exists some coupling $\nu_{dom}$
between $\nu_1=\pi_1^*\left(\nu_{dom}\right)$ and $\nu_2=\pi_2^*\left(\nu_{dom}\right)$
which is supported by $\Omega_\leq := \left\{\left(b_1,b_2\right)
\in \Omega_{\{0,1\}}^2,b_1\leq b_2\right\}$.

We can extend this notion to classes of measures.
Let $g:[0,1]\to [0,1]$ be a non-decreasing function,
right-continuous in $0$ with $g(0)=0$.
We say that the class $\mathcal{N}$ is $g$-dominated by $\mathcal{N}'$ if,
for any $\epsilon>0$ and any $\epsilon$-noise $\nu\in\mathcal{N}$,
there exists a $g(\epsilon)$-noise $\nu'\in\mathcal{N}'$ such that $\nu'\geq\nu$.
\end{definition}

Using this domination property,
we can then prove the following result,
that most notably does not depend on the distance $d$ used for the stability.

\begin{proposition} \label{prop:stableDomination}

If the SFT $\Omega_\F$ is $f$-stable on the class $\mathcal{N}'$ for the distance $d$,
and $\mathcal{N}$ is $g$-dominated by $\mathcal{N}'$,
then $\Omega_\F$ is $(f\circ g)$-stable on the class $\mathcal{N}$ for the distance $d$.

% Here an alternate proof :
%
% \begin{proof}
% Consider a measure $\lambda\in \widetilde{\M_\F^{\mathcal{N}}}(\epsilon)$,
% a coupling between $\mu\in\pi_1^*(\lambda)$ and
% the $\epsilon$-noise $\nu=\pi_2^*(\lambda)\in\mathcal{N}$.
% As $\mathcal{N}$ is $g$-dominated by $\mathcal{N}'$,
% there exists a coupling $\lambda''$ between $\nu$
% and a $g(\epsilon)$-noise $\nu'\in\mathcal{N}'$ such that $\nu'\geq\nu$.
% 
% As in the proof of Lemma~\ref{lem:BesicovitchDistance},
% we can here ~\cite[Chapter 6]{Gla03} couple
% the couplings $\lambda$ and $\lambda''$ that are ``compatible''
% in the sense that they share the projection $\nu$.
% Thus, there exists $\lambda'$ such that
% $\left(\pi_1,\pi_2\right)^*\left(\lambda'\right)=\lambda$
% and $\left(\pi_2,\pi_3\right)^*\left(\lambda'\right)=\lambda''$.
% 
% Consider $\left(\omega,b,b'\right)$ in the support of $\lambda'$.
% As $(\omega,b)$ is in the support of $\lambda$, it is an (obscured)
% locally admissible configuration of $\Omega_{\widetilde{\F}}$.
% As $\left(b,b'\right)$ is in the support of $\lambda''$,
% we have $b\leq b'$.
% Thus, $\left(\omega,b'\right)\in\Omega_{\widetilde{\F}}$
% must be locally admissible too.
% 
% In other words, we have $\left(\pi_1,\pi_3\right)^*\left(\lambda'\right)
% \in\widetilde{\M_\F^{\mathcal{N}'}}(g(\epsilon))$.
% We conclude the proof by using the $f$-stability on $\mathcal{N}'$,
% so that $d\left(\mu,\M_\F\right) \leq f(g(\epsilon))$.
% \end{proof}

\begin{proof}
Let us assume that $\Omega_\F$ is $f$-stable on the class $\mathcal{N}'$.
Consider $\mu\in\M_\F^\mathcal{N}(\epsilon)$.
If we prove that $\mu\in\M_\F^{\mathcal{N}'}(g(\epsilon))$,
then $d_B\left(\mu,\M_\F\right) \leq f(g(\epsilon))$ by $f$-stability.

In order to prove this,
let us consider a measure $\lambda\in\widetilde{\M_\F^{\mathcal{N}}}(\epsilon)$
such that $\pi_1^*(\lambda)=\mu$,
with $\pi_2^*(\lambda)=\nu\in \mathcal{N}$ an $\epsilon$-noise.
By domination, there exists a $g(\epsilon)$-noise $\nu'\in\mathcal{N}'$
such that $\nu'\geq\nu$, with a coupling $\nu_{dom}$ between them.

Then, using the disintegration theorem, for $\nu$-almost-any $b\in\Omega_{\{0,1\}}$,
there is a measure $\mu_b$ on $\Omega_\A$ such that,
for any two Borel sets $A\subset\Omega_\A$ and $B\subset \Omega_{\{0,1\}}$:
\[
\lambda\left(A\times B\right)=\int_B \mu_b(A)\mathrm{d}\nu(b)
=\int \mu_b(A) \mathds{1}_B(b) \mathrm{d}\nu_{dom}\left(b,b'\right) .
\]
Now, we can naturally define the measure $\lambda'$
on $\Omega_\A\times\Omega_{\{0,1\}}$ as :
\[
\lambda'\left(A\times B\right)=
\int \mu_b(A) \mathds{1}_B\left(b'\right) \mathrm{d}\nu_{dom}\left(b,b'\right) .
\]
By taking $B=\Omega_{\{0,1\}}$,
it is clear that $\pi_1^*\left(\lambda'\right)=\pi_1^*(\lambda)=\mu$.
Now, by taking $A=\Omega_\A$,
we conclude that $\pi_2^*\left(\lambda'\right)=\pi_2^*\left(\nu_{dom}\right)=\nu'$.
Moreover, consider
$\widetilde{w}=\left(w,0^{I(w)}\right)\in\widetilde{\F}$ a forbidden pattern.
Since $\nu_{dom}$ is supported
by $\Omega_{\leq}=\left\{\left(b,b'\right), b\leq b'\right\}$:
\[
\lambda'\left(\left[\widetilde{w}\right]\right) =
\int \mu_b([w]) \mathds{1}_{0^{I(w)}}\left(b'\right) \mathrm{d}\nu_{dom}\left(b,b'\right)
\leq \int \mu_b([w]) \mathds{1}_{0^{I(w)}}(b) \mathrm{d}\nu_{dom}\left(b,b'\right)
=\lambda \left(\left[\widetilde{w}\right]\right)=0 .
\]
Thence, $\lambda'$ is supported by $\Omega_{\widetilde{\F}}$,
so that $\lambda'\in\widetilde{\M_\F^{\mathcal{N}'}}(g(\epsilon))$.
At last, we demonstrated
that $\mu=\pi_1^*\left(\lambda'\right)\in\M_\F^{\mathcal{N}'}(g(\epsilon))$,
which concludes the proof.
\end{proof}
\end{proposition}

Now, in order to use this result,
we need to dominate $\gamma_n(\B)\subset \mathcal{D}_{2n+1}$ by $\B$.
By adapting a classical result, we can obtain the following bound:
\begin{proposition}[{{\cite[Theorem 1.3]{LiSchoSta97}}}]
The  $k$-dependent noise class $\mathcal{D}_k$
is polynomially $g_k$-dominated by $\B$,
with $g_k(\epsilon)\leq C\epsilon^{1/(2k+1)^d}$ for some constant $C$
that does \emph{not} depend on $k$ nor $d$.
\end{proposition}

\begin{corollary}
If the SFT $\Omega_\F$ is stable (resp. polynomially stable) on the Bernoulli class $\B$,
then it is also stable (resp. polynomially stable) on any dependent class $\mathcal{D}_k$.

Under the further assumption that there
exists a conjugacy $\theta:\Omega_\F\to\Omega_\G$,
then $\Omega_\G$ is also stable (resp. polynomially stable) on the class $\B$.

\begin{proof}
For the first part of the result, assume that $\Omega_\F$ is $f$-stable on $\B$.
As the class $\mathcal{D}_k$ is $g_k$-dominated by $\B$,
we may apply Proposition~\ref{prop:stableDomination},
so the SFT is $\left(f\circ g_k\right)$-stable on $\mathcal{D}_k$.
In particular, for the polynomial case,
if $f$ is a $O\left(\epsilon^\alpha\right)$,
then $f\circ g_k$ is a $O\left(\epsilon^{\alpha/(2k+1)^d}\right)$, still polynomial.
For the second part of the result,
we may use Theorem~\ref{thm:Conjugacy},
as $\Omega_\F$ is now $\left(f\circ g_{2r_{\theta^{-1}}+1}\right)$-stable
on $\gamma_{r_{\theta^{-1}}}(\B)\subset \mathcal{D}_{2r_{\theta^{-1}}+1}$.
In particular,
if $f\circ g_{2r_{\theta^{-1}}+1}$ is
a $O\left(\epsilon^{\alpha/\left(4r_{\theta^{-1}}+3\right)^d}\right)$,
then so is $D_\theta f\circ g_{2r_{\theta^{-1}}+1}
\left(C_{\theta^{-1}}\epsilon\right)+E\epsilon$.
\end{proof}
\end{corollary}

Notice how, because of the domination,
we are unable to preserve \emph{linear} stability.
Still, we have proven that stability on the class $\B$
is a conjugacy-invariant property.
In particular, stability on the class $\B$ is an intrinsic property
of a SFT $\Omega_\F$, which actually does not depend on the set of forbidden patterns $\F$
used to describe it.

As $g_k(\epsilon)\approx\epsilon^{1/(2k+1)^d}\underset{k\to\infty}{\longrightarrow} 1$
for any fixed value of $\epsilon$,
we conclude that even though a SFT stable on $\B$
is stable on all the classes $\mathcal{D}_k$,
this stability does not reach the limit class $\mathcal{D}=\bigcup_{k\in\N}\mathcal{D}_k$,
which would be the most natural generalisation of $\B$ stable under
all the morphisms $\gamma_k$.

\begin{remark}[Other Classes of Noise]
So far, we have only talked about noises in the class $\mathcal{D}$
of finite-range dependence, which we brought back to the independent case $\B$.
This focus is purposeful,
as pretty much all our further stability results will be proven on the class $\B$.

If we consider infinite-range dependencies, then we allow periodic noises,
\emph{i.e.} noises defined as uniform laws among the finite set of translations
of a periodic configuration $b\in\Omega_{\{0,1\}}$.
In most of the interesting cases, the rigid structure of such noises allows us to
explicitly construct measures that do not converge to $\M_\F$ for $d_B$,
as in Subsections~\ref{subsec:UnstablePeriod} and~\ref{subsec:GridNoise}.

The remaining in-between case would be that of infinite-range dependencies
but with correlations that decrease and go to $0$ as the distance goes to $\infty$.
This case notably encompasses the Gaussian Free Field,
as well as some Gibbs measures.
This may be the most physically realistic case,
but is also the harder to study,
so we will set it aside for the rest of this exploratory work.
\end{remark}

\section{Classification of the One-Dimensional Stability} \label{sec:1D}

Now that we have proved a general conjugacy invariance of the stability,
let us focus on a more specific framework, the one-dimensional (1D) case.
This case has already been widely studied,
and since the set of configurations of a 1D SFT can be seen
as the set of bi-infinite paths in an word automaton,
a lot of properties have been classified~\cite{LiMa21}. 

The section will be concluded by a discussion
on how to transpose general SFTs from $d$ to $d+1$ dimensions
while preserving their (un)stable behaviour;
this subsection is more technical and may be skipped without
harming the reading of the rest of the article.

We will now briefly introduce the main tool allowing for such a classification,
word automata, and then use it to classify stability
as a consequence of the aperiodicity of an automaton.
To put it shortly, stability of the SFT will be roughly equivalent to
the uniqueness of a communication class in the automaton, which must be aperiodic.

\subsection{1D SFTs and Word Automata}

In the 1D case, patterns and configurations are also called \emph{words}.
Because of their linear structure, words exhibit some automatic properties
not encountered in higher dimensions.

\begin{definition}[Diameter of a Set of Forbidden Patterns]
For a window of cells $I\subset \Z$, we denote $d(I)=\max(I)-\min(I)$ its diameter.
For a word $w\in \A^I$, $d(w)=d(I)$.
Finally, for a set of forbidden patterns $\F$,
we denote $d(\F)=\max_{w\in\F} d(w)$ its maximal diameter.
\end{definition}

Consider an automaton $G_\A^d$ where states are words in $\A^d$,
with transitions $au\overset{b}{\longrightarrow}ub$
for any $u\in \A^{d-1}$ and $a,b\in \A$.
Then it is equivalent to consider bi-infinite words $w\in\A^\Z$
and bi-infinite sequences of transitions in this word automaton.

Please note that this definition looks at words left-to-right,
but we could likewise look at right-to-left $ub\overset{a}{\longrightarrow}au$
transitions without changing any of the following (a)periodicity properties
nor the (in)stability results they imply.

\begin{definition}[Word Automaton]
Consider a set of forbidden words $\F$.
We define the automaton $G_\F$ induced
by restricting $G_\A^{d(\F)}$ to the states $w\in \A^d$
that contain no forbidden pattern,
that are locally admissible.
\end{definition}

Note that a configuration $w\in\Omega_\A$ corresponds to
a bi-infinite sequence of transitions of the automaton $G_\F$
if and only if $w\in\Omega_\F$ is a configuration of the SFT.

A SFT $\Omega_\F$ can be directly described by an automaton
instead of a set of forbidden patterns,
but we cannot \emph{construct} $G_\F$ out of $\Omega_\F$,
just as we cannot construct $\F$ itself.

As the number of states is finite, an infinite path exists
if and only if $G_\F$ contains a cycle, which allows us to decide
whether $\Omega_\F=\emptyset$ is empty or not in polynomial time.

\begin{definition}[Irreducible Automaton]
Two states $u,v\in \A^d$ of $G_\F$ communicate if there is a path
from $u$ to $v$ and $v$ to $u$ in the directed graph induced by $G_\F$.

This gives us a partial equivalence relation, whose classes
are the communication classes.
As long as $\Omega_\F\neq \emptyset$,
there is a cycle in $G_\F$ so such a class always exists.

We say that $G_\F$ is irreducible if this class is unique.
Please note that this does not imply that \emph{all} the states of $G_\F$
are in the class.
For example, in the directed graph represented by
$a\rightarrow b\;\rotatebox[origin=c]{90}{$\circlearrowright$}$,
$\left\{b\right\}$ is the only communication class,
because there is no path from $b$ to $a$.
\end{definition}

\begin{definition}[Periodic Automaton]
Consider $G_\F$ an irreducible automaton.
We say that it is $p$-periodic if $p$ is the greatest common divisor
of the lengths of all the cycles found inside $G_\F$.
$G_\F$ is \emph{aperiodic} if $p=1$.

In the $p$-periodic case, there exists a partition
$C=\bigsqcup_{j\in\Z/p\Z} C_j$ of the communication class
such that for any transition $u\to v$ of the automaton we must have
$u\in C_j$ and $v\in C_{j+1}$ for some $j\in\Z/p\Z$.
\end{definition}

\subsection{A Uniquely Ergodic Unstable Example} \label{subsec:UnstablePeriod}

For a stable SFT, as $\epsilon$ goes to $0$,
a generic noisy configuration has arbitrarily few differences
with a generic clear configuration of the SFT.
In the specific case of uniquely ergodic SFTs,
since there is only one measure in $\M_\F$,
a stronger structure is expected for generic clear configurations,
hence a prior motivation to study this case in particular.

In the 1D case, uniquely ergodic SFTs
are reduced to the finite orbit of a periodic configuration.
Hence, consider the simplest non-trivial uniquely ergodic 1D SFT,
whose only two configurations are $\omega_0=(01)^\Z$ and $\omega_1=(10)^\Z$
(such that $\omega_i(k)\equiv k+i[2]$).
This system is induced by the forbidden patterns $\F=\{00,11\}$,
it admits a unique invariant measure (hence it is uniquely ergodic),
and it is irreducible $2$-periodic.

We define the $p$-periodic noise $\nu_p$
whose configurations are the $p$ translations of $\left(0^{p-1}1\right)^\Z$.
With this noise, $\nu_p([1])=\frac{1}{p}$ goes to $0$ as $p\to\infty$.

Consider then $\lambda_p \in \widetilde{\M_\F}\left(\frac{1}{p}\right)$
such that $\pi_2^*\left(\lambda_p\right)=\nu_p$, and on each clear window of size $p-1$,
we use alternatively the restriction of $\omega_0$ or $\omega_1$.
Up to the values under obscured cells,
which will bear no influence on the following proposition,
we may assume without loss of generality
that $\lambda_p$ is supported by $2p$-periodic configurations.

\begin{proposition}
We have $d_B\left( \pi_1^*\left(\lambda_p\right),\M_\F\right)=
\frac{1}{2}-O\left(\frac{1}{p}\right)$.

\begin{proof}
Consider $(w,b)$ a $2p$-periodic configuration for $\lambda_p$,
and an interval $I=\left\llbracket k,k+2p\right\llbracket$ of size $2p$.
The restriction of $w$ in the window must coincide with the restriction of
$\omega_0$ in (at least) $p-1$ cells,
and \emph{cannot} coincide with $\omega_0$
on the $p-1$ cells specifically aligned with $\omega_1$,
thus $d_{2p}\left(w|_I,\omega_0|_I\right) \geq \frac{p-1}{2p}
= \frac{1}{2}-\frac{1}{2p}$.
More generally, for any choice of $n=2pq+r$ with $0\leq r<2p$,
and any interval $I$ of size $n$, which contains $q$ distinct intervals of size $2p$,
$d_n \left(w|_I,\omega_0|_I\right) \geq \frac{q(p-1)}{n}=\frac{p-1}{2p+\frac{r}{q}}
\underset{n\to\infty}{\longrightarrow}\frac{1}{2}-\frac{1}{2p}$,
which naturally gives the limit bound 
$d_H\left(w,\omega_i\right)\geq \frac{1}{2}-\frac{1}{2p}$.

As the lower bound on $d_H$ holds for
any configuration $(w,b)$ in the support of $\lambda_p$
and both globally admissible configurations $\omega_0$ and $\omega_1$,
it extends to $d_B\left(\pi_1^*\left(\lambda_p\right),\M_\F\right)$.
\end{proof}
\end{proposition}

We can generalise this result to all periodic SFTs without much effort,
provided we use periodic noises of the form
$\left(0^p 1^d\right)^\Z$ with $d\geq d(\F)$ and $p\to\infty$.
We will instead exhibit another $d_B$-instability,
but with Bernoulli noises in a further subsection.

\subsection{Irreducible Aperiodic Stability}

In a 1D setup, as long as $\Omega_F\neq \emptyset$,
there is always a cycle in the word automaton, thus a periodic configuration.
The aperiodicity of the automaton only implies the
\emph{existence} of aperiodic configurations in $\Omega_\F$,
which will prove to be sufficient to obtain stability.

Let us denote $L\left(\Omega_\F\right)$ the language of the SFT,
the set of words in $\A^*$ that are a restriction of a configuration of $\Omega_\F$.

\begin{remark} \label{rmk:StickingWords}

Consider a set of forbidden words $\F$ such that the automaton $G_\F$
has a unique communication class which is aperiodic.
It easily follows from the aperiodicity of $G_\F$ that
there exists a constant $n_0\in \N$ such that, for any $u,v\in L\left(\Omega_\F\right)$
and $n\geq n_0$, there exists a word $w\in \A^n$
such that $uwv\in L\left(\Omega_\F\right)$.
This constant can easily be computed from $G_\F$ in polynomial time --
with respect to the size $\left( |\A|+\sum_{w\in \F} |w| \right)\in \N$ for example.

For more details on the basic properties of the 1D case,
one may refer to the classic book by Lind and Marcus~\cite{LiMa21}.
\end{remark}

Assuming we can cut an obscured configuration from $\Omega_{\widetilde{\F}}$
into globally admissible words all distant by at least $n_0$,
then we will be able to rewrite these gaps in order to
obtain a globally admissible configuration.
If we only exclude obscured cells, then we will obtain a sequence of
locally admissible clear words instead,
that may not be globally admissible,
and the gaps between these words may be too small.
By leaving out the $\left\lceil\frac{n_0}{2}\right\rceil$-neighbourhood
around each obscured cell,
we make sure the gaps are big enough to be fillable.

The following proposition is a stronger 1D version
of the reconstruction function $\phi$
described in Remark~\ref{rmk:LocalGlobalGeneric}.

\begin{proposition} \label{prop:1DReconstruction}

There exists a constant $C(\F)$ such that, for any locally admissible word $u\in \A^*$,
by removing (at most) $C$ letters on each end,
we obtain instead a globally admissible word $v\in L\left(\Omega_\F\right)$.

\begin{proof}
Note that a path of length $n$ in $G_\F$ visits $n+1$ windows,
and represents a word of length $d(\F)+n$.
Thus, we may assume that $C\geq \frac{d(\F)}{2}$,
so that we only need to consider words long enough to represent a finite path
in the automaton $G_\F$.

As long as we visit vertices in the communication class of $G_\F$,
we can infinitely extend the path on both directions,
thence the word we encode is globally admissible.

Issues arise when we visit other states, which explicitly correspond
to windows that never occur in a bi-infinite path, thus non-globally admissible words.
As there is only one communication class,
no path can cycle through such a state.
Hence, if there are $k$ states of $G_\F$
outside of the communication class,
by removing $k$ states on each end of the path,
we make sure that the path only visits the communication class,
thus corresponds to a globally admissible word.
Hence, $C=\max\left(k,\left\lceil\frac{d(\F)}{2}\right\rceil\right)$
is big-enough.

If we want a better constant,
we can replace $k$ by the maximum of
the length of the longest path among vertices
\emph{outside of} yet \emph{connected to} the communication class,
and half of the longest path not connected to the class.
\end{proof}
\end{proposition}

Just like $n_0$, $C$ can be computed from $G_\F$ in polynomial time.
Now, if we remove a $C$-neighbourhood around each obscured cell,
then we obtain a sequence of globally admissible words.
Finally, by removing a $D$-neighbourhood
with $D=\max\left(C,\left\lceil \frac{n_0}{2}\right\rceil\right)$,
we make sure that we obtain alternately globally admissible words
and fillable gaps.
This really is the key idea of the following theorem,
whose proof mostly aims at properly explaining why
the transformation we perform is a $\sigma$-invariant morphism that
returns a clear globally admissible configuration.

\begin{theorem} \label{thm:Characterisation1D}

A 1D SFT $\Omega_\F$ with an aperiodic automaton $G_\F$ is linearly stable.

\begin{proof}
In order to obtain linear stability, we will consider a measure
$\lambda\in\widetilde{\M_\F}(\epsilon)$, and build
a measurable mapping $\psi:\Omega_{\widetilde{\F}}\to\Omega_\F$,
so that $d_H\left(\pi_1(\omega),\psi(\omega)\right)$ is small.
Let us notice that $\psi$ does not need to be defined on $\Omega_{\widetilde{\F}}$,
but only on a high-probability support $S\subset \Omega_{\widetilde{\F}}$.
In such a case, we may add a third \emph{independent} coordinate to $\lambda$
that follows some given law in $\M_\F$, and project onto this coordinate with $\psi$
outside of the event $S$. This way, we have:
\[
d_B\left(\pi_1^*(\lambda),\M_\F\right) \leq 
\int\limits_S d_H\left(\pi_1(\omega),\psi(\omega)\right) \mathrm{d}\lambda(\omega)
+ \lambda\left(S^c\right) .
\]

Consider the cellular automaton $\gamma_D$ on $\Omega_{\{0,1\}}$,
as defined in Definition~\ref{def:thickNoise}.
This morphism obscures the cells in the $D$-neighbourhood,
as described in Definition~\ref{def:thickNoise}.
This process is clearly measurable,
and naturally extends as a morphism on $\Omega_{\widetilde{\F}}$.
We now need to map this subset of
$\Omega_{\widetilde{\F}}$ into $\Omega_\F$ in a measurable way.

The issue now is that, while we can manually fill each gap,
issues may arise with the order of the operations.
Indeed, assume we decide on a word $w\left(u_1,u_2,n\right)$
for any words $u_1,u_2\in L\left(\Omega_\F\right)$ and any gap
of size $n\geq n_0$, as in Remark~\ref{rmk:StickingWords}.
Naturally, if we have \emph{three} globally admissible words $u_1$,
$u_2$ and $u_3$ as well as two gaps $i$ and $j$,
then we can fill the leftmost gap first,
with $v=u_1 w\left(u_1,u_2,i\right)u_2$ and then the second one
with $v w\left(v,u_3,j\right) u_3$.

There are several ways to proceed, but we chose here
to be able to fill those gaps simultaneously,
so that the $\sigma$-invariance of the morphism directly follows.
To ensure we can fill the gaps simultaneously,
we simply need to know the leftmost and rightmost states
of $G_\F$ corresponding to $u_2$,
which requires in turn $\left|u_2\right| \geq d(\F)$.
By looking at the $\left\lceil \frac{d(\F)}{2}\right\rceil$-neighbourhood of a clear cell
in a configuration of $\gamma_D\left(\Omega_{\{0,1\}}\right)$,
we can see whether it belongs to a long-enough globally admissible clear word,
and obscure it if it does not.
Let us name $\theta$ the cellular automaton on $\Omega_{\{0,1\}}$
obtained by applying $\gamma_D$
and then this new measurable process.
We identify $\theta$ with the morphism on $\Omega_{\widetilde{\F}}$
that leaves the first coordinate unchanged.

For a configuration $\omega\in\Omega_{\widetilde{\F}}$,
the obscured cells in $\theta(\omega)$ are all in a $E$-neighbourhood
(with $E=D+\left\lceil\frac{d(\F)}{2}\right\rceil$)
of the original obscured cells,
so we still have a linear control on the frequency of obscured cells.

Now, all clear words of an obscured configuration $\theta(\omega)$
are of length at least $d(\F)$.
For such words, we can define $w\left(u_1,u_2,n\right)$
using only the $d$ rightmost letters of $u_1$ and the $d$ leftmost letters of $u_2$,
which won't change if we change letters on the other end of $u_1$ or $u_2$.
Thus, we can simultaneously replace all the obscured windows
by the corresponding clear words.
Let us name $\psi(\omega)$ the configuration obtained now.

There is one last issue to deal with,
\emph{i.e.}\ the fact that $\psi(\omega)$ consists of one big globally admissible
clear word, but that it may have an infinite obscured window on the left or the right.
Let us name $S$ the set of configurations where this phenomenon does not happen.
So far, we obtained a set $S$ and defined a morphism $\psi:S\to \Omega_\F$,
as stated in the first paragraph of the proof,
so let us now study the two terms of the bound.

First, inside of $S$, $d_H\left(\pi_1(\omega),\psi(\omega)\right) \leq 
d_H\left(\pi_2(\theta(\omega)),0^\infty\right)\leq 
(2E+1)d_H\left(\pi_2(\omega),0^\infty\right)$.
Thus, $\int_S d_H\left(\pi_1(\omega),\psi(\omega)\right)\mathrm{d}\lambda(\omega)
\leq (2E+1)\int d_H\left(\pi_2(\omega),0^\infty\right)\mathrm{d}\lambda(\omega)$.
So far, the bound holds for any configuration in $\Omega_{\widetilde{\F}}$,
any measure $\lambda\in\M_{\widetilde{\F}}$.

Assume now that $\lambda\in\widetilde{\M_\F}(\epsilon)$.
Then, using Birkhoff's pointwise ergodic theorem,
$\int d_H\left(\pi_2(\omega),0^\infty\right)\mathrm{d}\lambda(\omega)
= \pi_2^*(\lambda)([1])\leq \epsilon$.
Now we just need to study $\lambda\left(S^c\right)$ to conclude.
By symmetry, up to an added factor $2$,
$\lambda\left(S^c\right) \leq 2 \lambda\left(T\right)$ where
$T$ is the event where there is a infinite obscured window on the right
in the configuration $\theta(\omega)$.
If $\omega\in T$, then $1$ must at least have a $\frac{1}{2D+1}$ density in
the configuration $\pi_2(\omega)$ to begin with,
so that $\lambda(T)\times \frac{1}{2D+1} + \lambda\left(T^c\right)\times 0 
\leq \epsilon$ and $\lambda(T)\leq (2D+1)\epsilon$.

At last, we obtain the explicit bound
$d_B\left(\pi_1^*(\lambda),\M_\F\right) \leq 3(2E+1)\epsilon$,
with $E$ an explicit constant, computable in polynomial time.
\end{proof}
\end{theorem}

Remark that, when using independent $\epsilon$-Bernoulli noises,
as $\lambda(S)=1$, then we lose the factor $3$ in this upper bound,
but the constant is still in the same general order of magnitude.

\subsection{Periodic Instability}

In the previous subsection, we proved stability in
the aperiodic case.
The proof made full use of aperiodicity,
in the sense that the obscured cells can induce
a gap of arbitrary size into any globally admissible configuration,
and aperiodicity is needed to guarantee such a gap is fillable.
We also clearly saw how this approach failed
in the introductory $2$-periodic example,
using a periodic noise to precisely quantify
the amount of differences on finite windows
in order to obtain $d_B$-instability.
Our objective is now to prove a broader periodic $d_B$-instability result,
but for seemingly more natural noisy configurations,
using $\epsilon$-Bernoulli noises.

\begin{theorem}
Consider a SFT $\Omega_\F$ such that $G_\F$
is irreducible $p$-periodic ($p\geq 2$).
Then for any $\epsilon>0$
there exists $\mu_\epsilon\in \M_\F^\B(\epsilon)$
such that $d_B\left(\mu_\epsilon,\M_\F\right)\geq
\frac{p-1}{p d(\F)}-\epsilon$.

\begin{proof}
Let us begin by considering the partition of $\Omega_\F$
into $p$ sets $\left(\Omega_j\right)_{j\in \Z/p\Z}$ induced by the states of $G_\F$,
so that if $\omega\in\Omega_i$, then $\sigma_k(\omega)\in \Omega_{i+k}$.

Consider also once and for all a periodic word $\omega_0\in\Omega_\F$,
that corresponds to an infinite cycle of $G_\F$.
Note that this cycle may not be of length $p$ but a multiple of it --
\emph{e.g.}\ if $G_\F$ is made of a $6$-cycle and a $10$-cycle joined in a vertex,
it is $2$-periodic but has no $2$-cycle.
What matters is that $\omega_0$ has a finite orbit under translations.
What is more, by looking at a window of size $d(\F)$ of a translation of $\omega_0$,
we can identify to which state of $G_\F$ it corresponds and thus
deduce to which class $\Omega_j$ the translated configuration belongs to.
To construct $\mu_\epsilon\in\M_\F^\B(\epsilon)$,
consider the measure $\lambda_\epsilon$ obtained by:
\begin{enumerate}
\item taking the independent Bernoulli noise $\B(\epsilon)^{\otimes\Z}$ first,
\item identifying intervals of consecutive obscured cells, of length at least $d(\F)$,
and writing down letters of $\A$ uniformly at random under each such block,
\item in-between two such intervals, in a window that must have a clear cell on each end
and may contain some short obscured blocks in the middle,
we choose uniformly at random a translation of $\omega_0$ to write it down on the cells,
whether clear or obscured.
\end{enumerate}
It is apparent that this measure has an $\epsilon$-Bernoulli noise,
and that it is $\sigma$-invariant by construction.
The measure $\lambda_\epsilon$ is also strongly mixing, thus ergodic.
Indeed, consider two finite windows $I,J\subset \Z$
such that $\min (J) -\max(I) = n> d(\F)$.
Conditionally to the fact that the window $\llbracket \max(I)+1,\min(J)-1\rrbracket$
contains an obscured window of size $d(\F)$,
the windows $I$ and $J$ behave independently from each other.
As the probability of having such an obscured window goes to $1$ as $n\to\infty$,
we deduce the mixing property on cylinders,
so that $\lambda_\epsilon$ itself is strongly mixing.

Hence, if we cut down $\Z$ into consecutive windows of length $d(\F)$,
we obtain a measure on $\left( \A^d\times \{0,1\}^d\right)^\Z$.
This induced measure is \emph{also} $\sigma$-invariant and strongly mixing (thus ergodic),
so that Birkhoff's pointwise ergodic theorem applies.
Hence, the frequency of a $d(\F)$-interval in a configuration is
$\lambda_\epsilon$-almost-surely equal to its probability under $\lambda_\epsilon$.

A clear $d(\F)$-interval has probability $(1-\epsilon)^{d(\F)}$ of happening,
which we bound below by $1-d(\F)\epsilon$.
Under such an event, by construction, we can identify the state of $G_\F$ it represents,
thus to which class $\Omega_i$ it comes from.
Note that on such a clear window,
if $\omega$ and $\omega'$ belong to different classes $\Omega_i$ and $\Omega_j$,
then in particular they correspond to different states of $G_\F$ thus must differ in
\emph{at least} one cell.

Thus, for \emph{any} globally admissible state $\omega\in\Omega_\F$
and a \emph{typical} locally admissible state
$\left(\omega',b\right)\in\Omega_{\widetilde{\F}}$
under $\lambda_\epsilon$, we have:
\[
d_H\left(\omega,\omega'\right) \geq (1-d(\F)\epsilon)
\times \frac{p-1}{p}\times \frac{1}{d(\F)} .
\]
The first factor comes from the frequency of clear windows,
the second one from the probability of $\omega'$ not being in the same class as $\omega$
conditionally to some clear window,
and the third one from the minimal number of differences in such a window
of size $d(\F)$ under the previous event.

It immediately follows that
$d_B\left( \mu_\epsilon,\M_\F\right) \geq \frac{p-1}{pd(\F)}- \frac{p-1}{p}\epsilon$,
which concludes the proof.
\end{proof}
\end{theorem}

\begin{remark}
Note that this instability result
can be adapted from the periodic case to the non-irreducible case
where there are several communication classes,
by using finite trajectories evolving inside distinct communication
classes instead of words ``aligned'' along different periods of the system,
even if all the classes are aperiodic.
In such a framework, if $p$ denotes the number of classes,
then we can obtain the very same lower bound as before.
We will omit the proof of such an assertion for the sake of brevity,
as it offers no further insight on the topic.
\end{remark}

\subsection{Extension to Higher Dimensions}

This subsection dives deeper into the intricacies of couplings
from a measure theory viewpoint,
which offers a different insight on the objects we are working on,
but can also be skipped by an unfamiliar reader
as it is independent from everything that will follow.

Our goal here is to provide a simple way to extend noisy SFTs into higher dimensions.
First of all, let us quickly characterise $\sigma$-invariant couplings,
which will be useful for the main result of this subsection.

\begin{lemma} \label{lem:InvariantCoupling}
Let $\lambda$ be any probability measure on $\Omega_\A\times\Omega_\A$,
with $\mu=\pi_1^*(\lambda)$.
We can factorise $\lambda(A\times B)=\int_A \nu_x(B)\mathrm{d}\mu(x)$,
such that $x\mapsto \nu_x(B)$ is measurable for any cylinder $B$,
and that $B\mapsto \nu_x(B)$ is a probability measure
for $\mu$-almost-any $x\in \Omega_A$.
This result is known as the disintegration theorem.

Assume now that $\mu$ is $\sigma$-invariant.
Then $\lambda$ is also $\sigma$-invariant if and only if
the equality $\nu_x(B)=\nu_{\sigma_k(x)}\left(\sigma_k(B)\right)$ holds
for any $k\in\Z^d$, any measurable cylinder $B$, and $\mu$-almost-any $x\in\Omega_\A$.

\begin{proof}
Consider two cylinders $A$ and $B$, as well as $k\in \Z^d$.
As stated, we have the equality $\lambda(A\times B)=\int_A \nu_x(B)\mathrm{d}\mu(x)$.
Likewise, as $\mu$ itself is $\sigma$-invariant:
\[
\lambda\left(\sigma_k(A\times B)\right)
=\lambda\left(\sigma_k(A)\times\sigma_k(B)\right)
=\int\limits_{\sigma_k(A)} \nu_x\left(\sigma_k(B)\right)\mathrm{d}\mu(x)
=\int\limits_A \nu_{\sigma_k(y)}\left(\sigma_k(B)\right)\mathrm{d}\mu(y) .
\]

Now, the measure $\lambda$ is $\sigma$-invariant if and only if,
for any cylinder $B$ and $k\in \Z^d$,
we have $\lambda(A\times B)=\lambda\left(\sigma_k(A\times B)\right)$
for any cylinder $A$.
Using the integral expressions, $\int_A\nu_x(B)\mathrm{d}\mu(x)=
\int_A \nu_{\sigma_k(x)}\left(\sigma_k(B)\right)\mathrm{d}\mu(x)$.
It is equivalent for this equality to hold for any $A$
and for the functions to be $\mu$-almost-surely equal,
which concludes the proof.
\end{proof}
\end{lemma}

Note how the measures $\nu_x$ are not necessarily $\sigma$-invariant.
In particular, using the Dirac measures $\nu_x=\delta_x$ -- which are obviously
not $\sigma$-invariant -- gives us a \emph{diagonal} coupling between $\mu$ and itself,
such that $\pi_2^*(\lambda)=\mu$ too, which is $\sigma$-invariant.

Now, given a $d$-dimensional SFT $\Omega_\F$,
it is possible to extend $\F$ into $\F'$ in $d+1$ dimensions,
by replacing every forbidden pattern $w\in\A^{I(w)}$ on the window $I(w)\subset \Z^d$
by $w'\in \A^{I(w)\times\{0\}}$ with $I(w)\times\{0\}\subset \Z^{d+1}$.
This way, $\Omega_{\F'}=\left\{ \left(\omega_i\right)_{i\in\Z},\forall i \in \Z,
\omega_i\in\Omega_\F\right\}$.
In other words, each slice (with a fixed last coordinate)
represents a copy of the original SFT,
with no constraints on how to align the slices.
In particular, if $\mu \in \M_\F(\epsilon)$,
then by coupling all these layers independently,
we obtain $\mu^{\otimes \Z} \in \M_{\F'}(\epsilon)$.

Let us now prove that (in)stability of a SFT is in some sense preserved
through this transformation.
Thus, as we exhibited (un)stable 1D examples earlier in this section,
this will imply the existence of (un)stable systems in \emph{any} dimension.

Consider the projection $\zeta : b\in\{0,1\}^{\Z^{d+1}} \mapsto
b|_{\Z^d\times \{0\}} \in\{0,1\}^{\Z^{d}}$,
that commutes with translations in $\Z^d$.
More generally, we will use $\zeta$ as a multipurpose projector
for any alphabet $\A$ instead of $\{0,1\}$.
For a given class of $(d+1)$-dimensional noises $\mathcal{N}'$,
we obtain the $d$-dimensional class $\mathcal{N}=\zeta^*\left(\mathcal{N}'\right)$.
In particular, if $\mathcal{N}'$ is the class of $(d+1)$-dimensional Bernoulli noises,
then $\mathcal{N}$ is the class of $d$-dimensional Bernoulli noises.

To make things easier to read, we will distinguish the Besicovitch distances
$d_B$ in $d$ dimensions and $d_B'$ in $d+1$ dimensions (resp. $d_H$ and $d_H'$).

\begin{proposition}
Using the previous notations:
\begin{enumerate}
\item For any $\mu'\in\M_{\F'}^{\mathcal{N}'}(\epsilon)$,
we have $\mu=\zeta^*\left(\mu'\right)\in\M_\F^{\mathcal{N}}(\epsilon)$.
\item For any $\mu\in\M_\F^{\mathcal{N}}(\epsilon)$,
there exists $\mu'\in\M_{\F'}^{\mathcal{N}'}(\epsilon)$
such that $\mu=\zeta^*\left(\mu'\right)$.
\item In both cases, $d_B\left(\mu,\M_\F\right)=d_B'\left(\mu',\M_{\F'}\right)$.
\end{enumerate}

\begin{proof}
In all cases, going from dimension $d+1$ to dimension $d$ is a mere matter
of projection through $\zeta$, whereas going from dimension $d$ to $d+1$
is a bit trickier and will require us to make use of Lemma~\ref{lem:InvariantCoupling}.

First,
assume there is $\lambda_{noise}'\in\widetilde{\M_{\F'}^{\mathcal{N}'}}(\epsilon)$
such that $\mu'=\pi_1^*\left(\lambda_{noise}'\right)\in\M_\F$.
Thus, $\lambda_{noise}=\zeta^*\left(\lambda'_{noise}\right)
\in\widetilde{\M_\F^\mathcal{N}}(\epsilon)$,
so that:
\[
\mu=\zeta^*\left(\mu'\right)=
\zeta^*\left(\pi_1^*\left(\lambda_{noise}'\right)\right)
=\pi_1^*\left(\lambda_{noise}\right)\in\M_\F^\mathcal{N}(\epsilon) .
\]
This proves the first item.

Conversely, consider 
$\lambda_{noise}\in\widetilde{\M_\F^{\mathcal{N}}}(\epsilon)$ 
such that $\mu=\pi_1^*\left(\lambda_{noise}\right)$,
and let us build the desired measure $\mu'$.
Using Lemma~\ref{lem:InvariantCoupling},
we have $\mathrm{d}\lambda_{noise}(\omega,b)=\mathrm{d}\mu_b(\omega)\mathrm{d}\nu(b)$
with $\nu=\pi_2^*\left(\lambda_{noise}\right)\in\mathcal{N}$ an $\epsilon$-noise,
and $\mathrm{d}\mu_{\sigma_k(b)}\left(\sigma_k(\omega)\right)=\mathrm{d}\mu_b(\omega)$
for $\nu$-almost-any $b\in\Omega_{\{0,1\}}$.
As $\nu\in\mathcal{N}=\zeta^*\left(\mathcal{N}'\right)$,
there is $\nu'\in\mathcal{N}'$ such that $\nu=\zeta^*(\nu')$.
In particular, $\nu'$ is $\sigma$-invariant and must be an $\epsilon$-noise too.
Now, for any families of $d$-dimensional layers
$\omega'=\left(\omega_i\right)_{i\in\Z} \in \Omega_{\F'}$
and $b'=\left(b_i\right)$,
we define the measures $\mathrm{d}\mu_{b'}'\left(\omega'\right)=
\prod_{i\in\Z}\mathrm{d}\mu_{b_i}\left(\omega_i\right)$
and then $\mathrm{d}\lambda_{noise}'\left(\omega',b'\right)=
\mathrm{d}\mu_{b'}\left(\omega'\right)\mathrm{d}\nu'\left(b'\right)$.
Naturally, the measures $\mu_{b'}'$ are $\sigma_{e_{d+1}}$-invariant
-- invariant by translations on the last coordinate -- by construction,
and satisfy the criterion of Lemma~\ref{lem:InvariantCoupling}
because the measures $\mu_b$ did.
Thence, $\lambda_{noise}'$ is $\sigma$-inviariant,
so that $\lambda_{noise}'\in \widetilde{\M_{\F'}^{\mathcal{N}'}}(\epsilon)$.
At last, $\mu'=\pi_1^*\left(\lambda_{noise}'\right)$ is
such that $\zeta^*\left(\mu'\right)=\mu$, which proves the second item.

Finally, consider any two $\sigma$-invariant measures $\mu'$
and $\mu=\zeta^*(\mu)$.
We begin with the easier inequality,
by considering $\lambda'$ a coupling between $\mu'$ and $\nu'\in\M_{\F'}$
such that $d_B'\left(\mu',\M_{\F'}\right)=d_B'\left(\mu',\nu'\right)
=\int d_H'\left(x',y'\right)\mathrm{d}\lambda'\left(x',y'\right)$.
Using Birkhoff's pointwise ergodic theorem, this is equal to
$d_B'\left(\mu',\nu'\right)=\int\mathds{1}_{
x_0'\neq y_0'}\left(x',y'\right)\mathrm{d}\lambda'\left(x',y'\right)$.
Likewise, $\lambda=\zeta^*\left(\lambda'\right)$ is a coupling
between $\mu$ and $\nu=\zeta^*\left(\nu'\right)\in\M_\F$,
not necessarily such that $d_B$ is reached, but:
\[
d_B\left(\mu,\M_\F\right)\leq d_B\left(\mu,\nu\right) \leq
\int\mathds{1}_{x_0\neq y_0}\mathrm{d}\lambda(x,y)=
\int\mathds{1}_{x_0'\neq y_0'}\mathrm{d}\lambda'\left(x',y'\right)
=d_B'\left(\mu',\M_{\F'}\right) .
\]
For the reverse inequality, consider
$\lambda$ a coupling between $\mu$ and $\nu\in\M_{\F}$
such that $d_B\left(\mu,\M_\F\right)=\int\mathds{1}_{x_0\neq y_0}\mathrm{d}\lambda$.
As in Lemma~\ref{lem:InvariantCoupling},
we can factorise $\mathrm{d}\lambda(x,y)=\mathrm{d}\nu_x(y)\mathrm{d}\mu(x)$.
With $x',y'\in\Omega_{\F'}=\Omega_\F^\Z$, we define the family of measures
$\mathrm{d}\nu'_{x'}\left(y'\right)=\prod_{i\in\Z}\mathrm{d}\nu_{x_i}\left(y_i\right)$,
and then $\mathrm{d}\lambda'\left(x',y'\right)
=\mathrm{d}\nu'_{x'}\left(y'\right)\mathrm{d}\mu'\left(x'\right)$.
The measures $\nu'_{x'}$ satisfy the criterion of Lemma~\ref{lem:InvariantCoupling},
so that $\lambda'$ is $\sigma$-invariant.
This implies that it is a coupling between the measures $\mu'$
and $\nu'=\pi_2^*\left(\lambda'\right)\in\M_{\F'}$,
once again not necessarily optimal, such that:
\[
d_B'\left(\mu',\M_{\F'}\right)\leq d_B'\left(\mu',\nu'\right) \leq
\int\mathds{1}_{x_0'\neq y_0'}\mathrm{d}\lambda'\left(x',y'\right)=
\int\mathds{1}_{x_0\neq y_0}\mathrm{d}\lambda(x,y)
= d_B\left(\mu,\M_\F\right) .
\]
This concludes the proof of the last item.
\end{proof}
\end{proposition}

\begin{corollary} \label{cor:ExtensionDimension}
The $d$-dimensional SFT $\Omega_\F$ is $f$-stable (resp. unstable)
on the class $\mathcal{N}$ if and only if the $(d+1)$-dimensional SFT $\Omega_{\F'}$
is $f$-stable (resp. unstable) on the class $\mathcal{N}'$.

\begin{proof}
Going from (un)stability on $\Omega_{\F'}$ to $\Omega_\F$
is once again a simple matter of projecting measures with $\zeta$
so we won't insist further on these implications.

If $\Omega_\F$ is $f$-stable on $\mathcal{N}$,
and $\mu'\in\M_{\F'}^{\mathcal{N}'}(\epsilon)$,
then $\mu=\zeta^*\left(\mu'\right)\in\M_\F^{\mathcal{N}}(\epsilon)$
using Item~1 of the previous proposition,
so that $d_B'\left(\mu',\M_{\F'}\right)=d_B\left(\mu,\M_\F\right)\leq f(\epsilon)$
using Item~3.
Thus, $\Omega_{\F'}$ is $f$-stable.

Now, if $\Omega_\F$ is unstable,
we have a sequence of measures $\mu_n\in\M_\F\left(\epsilon_n\right)$
with $\epsilon_n\underset{n\to\infty}{\longrightarrow}0$
such that $\inf_{n\in\N} d_B\left(\mu_n,\M_\F\right) = d>0$.
Then, with the measures $\mu_n'\in\M_{\F'}^{\mathcal{N}'}\left(\epsilon_n\right)$
given by Item~2,
we conclude that $\inf_{n\in\N} d_B'\left(\mu_n',\M_{\F'}\right)\geq d$ too
with Item~3, thence $\Omega_{\F'}$ also is unstable.
\end{proof}
\end{corollary}

Using this corollary, we can in particular extend the (un)stable 1D SFTs
we exhibited earlier in order to obtain (un)stable SFTs in \emph{any} dimension.
Of course, these examples are not really satisfactory
and we will now strive for other higher-dimensional examples
in the following sections of this paper.

\begin{remark}
Another way to extend SFTs is the direct product.
If we consider two $d$-dimensional SFTs $\Omega_\F$ on the alphabet $\A$
and $\Omega_{\F'}$ on the alphabet $\A'$,
then we can build the SFT $\Omega_\F\times\Omega_{\F'}$ on the alphabet $\A\times\A'$.

Let us note that, with $\omega_1,\omega_2\in\Omega_\A$ and
$\omega_1',\omega_2'\in\Omega_{\A'}$, we have the inequalities:
\[
d_H\left(\omega_1,\omega_2\right)\leq
d_H\left(\left[\omega_1,\omega_2\right],\left[\omega_1',\omega_2'\right]\right)
\leq d_H\left(\omega_1,\omega_2\right)+d_H\left(\omega_1',\omega_2'\right) .
\]

Thence, if $\Omega_\F$ (resp. $\Omega_{\F'}$) is $f$-stable (resp. $f'$-stable)
on the same class $\mathcal{N}$,
then the product $\Omega_\F\times\Omega_{\F'}$
is $\left(f+f'\right)$-stable on the class $\mathcal{N}$.
If one of the SFTs is unstable,
then the product is unstable with the same lower bound.
\end{remark}

\section{Stability of Periodic SFTs} \label{sec:2DPeriodic}

In this section, we will explore the notion of stability
for higher-dimensional (\emph{2D+}) periodic SFTs.
Here, we really mean periodicity of the SFT and its configurations,
not of some associated structure like the word automaton of Section~\ref{sec:1D}.
First, we will show how to obtain instability using a grid noise,
like we did in Subsection~\ref{subsec:UnstablePeriod} for the 1D case.
We will then focus on Bernoulli noises and prove that,
using a percolation argument for 2D+,
we have linear stability in this framework.

There are several non-equivalent notions of periodicity in the 2D+ case.
We will in this case consider the strongest notion of periodicity,
\emph{i.e.}\ the existence of $\Z$-independent vectors $x_1,\dots, x_d\in\Z^d$,
such that \emph{any} configuration $\omega\in \Omega_\F$ is invariant
under any translation among those ($\sigma_{x_i}(\omega)=\omega$ for
any $1\leq i\leq d$).
% Z-libre <=> Q-libre <=> Q-base
% <=> atteint toutes les Q-droites <=> toutes les Z-droites.
Equivalently, we can always assume that these $d$ vectors
align with the $d$ axes of $\Z^d$, so that we can actually
simply repeat a base pattern defined on a hyper-rectangle
along those $d$ base directions.

Up to an added redundancy along some of those axes,
we may even go one step further and assume the base pattern is defined
on a hypercube whose edge-length is the smallest common multiple
of those of the hyper-rectangle.
This added hypothesis will worsen the constants obtained in the following proofs,
but will make notations a bit lighter as a trade-off.

\subsection{Instability for Grid Noises} \label{subsec:GridNoise}

\begin{definition}[Grid Noise]
Consider $k,n\in\N^*$ two positive integers.
We define the base pattern $b_{k,n}$
on the hypercube $\llbracket 0,k+n-1\rrbracket^d$
such that $b(x)=1$ \emph{iff} $\min_{1\leq i \leq n} x_i < k$.
We then identify $b_{k,n}$ to the configuration obtained by
extending this base pattern in all directions.
We finally define the $\sigma$-invariant noise:
\[
\nu_{k,n} = \frac{1}{(k+n)^d}\sum\limits_{x\in\llbracket 0,k+n-1\rrbracket^d}
\delta_{\sigma_x\left(b_{k,n}\right)} .
\]
The probability of an obscured cell in this noise is $\left(\frac{k}{k+n}\right)^d$.
\end{definition}

Assuming $k$ is greater than the maximal diameter of the forbidden patterns of $\F$,
then two distinct clear hypercubes (both translations of $\llbracket 0,n-1\rrbracket^d$)
are insulated from each other, and can be tiled independently,
as no forbidden pattern could have cells in both windows.
We will work under this assumption from now on.

\begin{proposition}
For any (non-trivial) periodic SFT,
there exists a constant $\delta(\F)>0$ such that, for any $\epsilon>0$,
there is a measure $\mu\in\M_\F(\epsilon)$ at distance at least $\delta$
from $\M_\F$.

\begin{proof}
Let us assume that $\F$ is a periodic SFT,
on the base hypercube $\llbracket 0,N-1\rrbracket^d$.
Then two distinct configurations $\omega\neq \omega' \in \Omega_\F$
differ on \emph{at least} one cell in any translation of the $N$-hypercube.

By monotonicity, we only need to prove it for arbitrarily small values of $\epsilon$.
We will prove the result for the noises $\nu_{k,nN}$ as $n\to\infty$,
for which the frequency of obscured cells is equal to
$\epsilon_n=\left(\frac{k}{k+nN}\right)^d$,
such that $\epsilon_n\underset{n\to\infty}{\longrightarrow}0$.

What we mean here by non-trivial is that there exists a non-constant
configuration $\omega_0\in\Omega_\F$
such that $\omega_0 \neq \sigma_{e_j}\left(\omega_0\right)$
for some $1\leq j\leq d$,
thus $\Omega_0 := \left\{ \sigma_x\left(\omega_0\right), x\in\Z^d\right\}$
has between $2$ and $N^d$ elements.

We define the noisy measure
$\lambda\in \widetilde{\M_\F}\left(\epsilon_n\right)$ as follows:
\begin{itemize}
\item first, pick a noise grid at random following the measure $\nu_{k,nN}$,
\item under any obscured cell pick a letter uniformly at random,
\item then, independently from the noise, and independently on each clear hypercube,
pick a configuration $\omega\in\Omega_0$ uniformly at random,
and finally restrict it to the corresponding hypercube.
\end{itemize}

Consider a configuration $\omega\in \Omega_\F$, and $\left(\omega',b\right)$
in the support of $\lambda$.
Almost-surely, in $\omega'\in\Omega_\A$,
a proportion $\frac{1}{\left|\Omega_0\right|}\geq \frac{1}{N^d}$
of the clear hypercubes from $b$ contains each translation of $\omega_0$.
Hence, in a proportion greater or equal to
$\frac{\left|\Omega_0\right|-1}{\left|\Omega_0\right|}\geq \frac{1}{2}$,
the configuration chosen for this hypercube is \emph{not} $\omega$.
For such a clear window, as a translation of $\llbracket 0, nN-1\rrbracket^d$,
contains $n^d$ distinct translations of $\llbracket 0, N-1\rrbracket^d$.
On each such sub-hypercube, $\omega$ and $\omega'$ differ on at least one cell.
Finally:
\[
d_H(\omega,\omega')\geq \frac{1}{2}\times\left(\frac{n}{k+nN}\right)^d .
\]
This inequality holds almost surely under $\lambda$ for \emph{any}
configuration $\omega\in\Omega_\F$, so for big enough values of $n\geq k$,
we obtain the lower bound
$d_B\left( \pi_1^*(\lambda),\M_\F\right) \geq \frac{1}{2(N+1)^d}$.
\end{proof}
\end{proposition}

\subsection{From Noisy SFTs to Percolations}

In the 1D case, under a Bernoulli noise,
having room for aperiodicity was what helped us
correct defects in the noisy configurations from $\Omega_{\widetilde{\F}}$
in order to couple them with globally admissible configurations in $\Omega_\F$,
while intrinsic periodicity of the SFT was precisely what prevented stability.
Yet, in the 2D+ case, we will see that periodicity \emph{helps} stability
as long as most of the clear cells are connected to each other in an induced percolation.

Once again, let us consider a variant of the reconstruction function
described in Remark~\ref{rmk:LocalGlobalGeneric}.
Here, $\phi_\F:\N\to \N$ is a non-decreasing function such that,
for any integer $n\in\N^*$, $\phi(n)\geq n$ and
whenever $\omega\in \A^{B_{\phi(n)}}$ is a locally admissible pattern,
its restriction $\omega|_{B_n}$ is globally admissible.
For the 1D case, we proved in Proposition~\ref{prop:1DReconstruction}
that this function can always be chosen as $\phi(n)=n+c$ for some $c\in\N$.
This property allowed us to convert a locally admissible configuration into
a globally admissible clear one, up to some ``peeling'' around obscured cells,
in the case of aperiodic word automata.
What we now want is to transpose this argument into the 2D+ case,
using purposely the redundancy induced by the periodicity.

\begin{lemma} \label{lem:LocalGlobal}

Consider a 2D+ periodic SFT $\Omega_\F$.
There exists a constant $c(\F)\in \N$ such that,
for any \emph{connected} cell window $I\subset \Z^d$,
if $w\in \A^{I+B_c}$ is locally admissible,
then $w|_I$ is globally admissible.

\begin{proof}
As the SFT is periodic, consider $N$ the size of a base hypercube like before.

Let us begin with the case where $I=\left\{e\right\}$ is made of a single cell.
Assuming a pattern $w$ on the window $e+B_{\left\lceil \frac{N}{2}\right\rceil}$
is \emph{globally} admissible,
then $w$ actually is the restriction of a configuration $\omega_e\in\Omega_\F$
that coincides with $w$ on the window,
and in particular $\omega_e|_I=w|_I$.
Thus, it is sufficient to consider
$c=\phi\left(\left\lceil \frac{N}{2}\right\rceil\right)$,
such that whenever $w$ is locally admissible on $e+B_c$,
it is globally admissible on $e+B_{\left\lceil \frac{N}{2}\right\rceil}$
so the previous paragraph applies.

More generally, consider any connected window of cells $I$,
and $w\in \A^{I+B_c}$ a locally admissible pattern.
For any cell $e\in I$, we can likewise obtain a configuration
$\omega_e\in\Omega_\F$ such that,
on the domain $e+B_{\left\lceil \frac{N}{2}\right\rceil}$,
the pattern $w$ and the configuration $\omega_e$ coincide.

Consider now two neighbouring cells $e,f\in I$.
As we left a bit of margin to begin with,
the intersection $\left(e+B_{\left\lceil \frac{N}{2}\right\rceil}\right)\cap
\left(f+B_{\left\lceil \frac{N}{2}\right\rceil}\right)$ contains a $N$-hypercube,
thus the same base pattern for both $\omega_e$ and $\omega_f$,
so that we actually have equality $\omega_e=\omega_f$.

As $I$ is connected, by induction,
the pattern $w|_I$ is actually a restriction of $\omega_e$,
hence globally admissible.
\end{proof}
\end{lemma}

For a noisy configuration $(\omega,b)\in\Omega_{\widetilde{\F}}$
to be close to a globally admissible one,
we need a high-density connected window $I$ such that all cells in $I+B_c$ are clear.
If such a window occurs with high probability,
then we will be able to control the distance of a noisy measure to $\M_\F$.
Notice that this behaviour can be characterised by looking solely at the noise $b$,
by studying a site percolation on $\Z^d$.
This is what we will do in the next subsection.

\subsection{Study of the Thickened Percolation}

We consider here the site percolation on $\Z^d$,
with configurations $b\in\Omega_{\{0,1\}}$.
In our framework, the \emph{open} cells will be the clear ones, with value $0$,
and the \emph{closed} ones will be the obscured ones, with value $1$.

As we want specific properties on the ``thickness'' of the infinite component
of this percolation (the connected window $I$ such that $I+B_c$ is open),
we will induce an auxiliary percolation.
If $\nu\in \M_{\left\{0,1\right\}}$ defines a random percolation,
then the percolation on $\Z^d$ induced by the measure $\gamma_n^*(\nu)$
will be called the $n$-thickened $\nu$-percolation,
with the cellular automaton $\gamma_n$ from Definition~\ref{def:thickNoise}.

In the article \emph{Density and Uniqueness in Percolation}
~\cite[Theorem 2]{BurKea89}, it is shown that
under a condition of \emph{finite energy} on the measure $\mu$, defined below,
the percolation almost-surely has \emph{at most} one infinite connected component.
This property holds true for any Bernoulli noise in particular.

\begin{definition}[Finite Energy]
Consider $w\in \A^I$ a finite pattern.
For a measurable set $B$, we define $B_w=\left\{
\omega\in \Omega_\A,\exists\omega'\in B, 
\omega|_{I^c}=\omega'|_{I^c},\omega|_I=w\right\}$ which is also measurable.

A measure $\mu$ has finite energy if,
for any finite pattern $w$ and any measurable set $B$,
we must have $\mu\left(B_w\right)>0$ whenever $\mu(B)>0$.
\end{definition}

Please note that thickened measures cannot have the free energy property.
Indeed, a consequence of free energy is that any cylinder has a positive measure.
However, for a $n$-thickened percolation, we cannot have three adjacent cells
with the pattern $010$ in a configuration $\gamma_n(b)$,
as the presence of a $1$ in a $n$-hypercube of $b$ implies its presence
in the left-translated or right-translated hypercube.
The result can nonetheless be effortlessly adapted to the case of thickened measures,
and we will sketch its proof here for completeness.

\begin{figure}
\centering
\includegraphics[scale=.6]{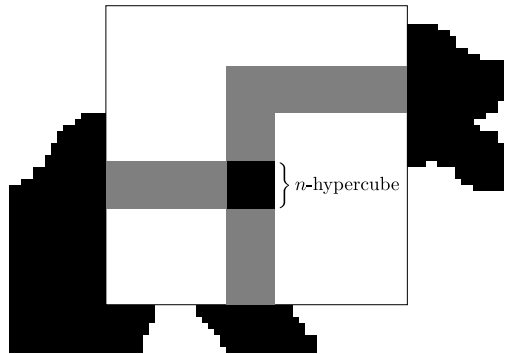}
\caption{Schematic representation of a trifurcation in the $n$-thickened case.}
\label{fi:Trifurcation}
\end{figure}

\begin{lemma}
When $\nu$ has the finite energy property, any thickened $\nu$-percolation
has \emph{at most} one infinite connected component.

\begin{proof}
The finite energy property still holds for the measures obtained through
the ergodic decomposition theorem, hence we can assume $\nu$ is ergodic.
As $\gamma_n$ is $\sigma$-invariant,
by definition of ergodicity, if $\nu$ is ergodic,
then so is the $n$-thickened $\nu$-percolation.

As a $\sigma$-invariant measurable function,
the number $N(b)$ of infinite components in the percolation $b$
is $\gamma_n^*(\nu)$-almost-surely constant.

If $N$ was infinite,
then for a big-enough hypercube $B$,
the probability of encountering three different infinite components in $\gamma_n(b)$
inside of it would be positive.

In the context of site percolation, a trifurcation of a configuration $b$ is
an open cell that is part of an infinite component, with exactly three open neighbours
such that if the cell was closed then these neighbours would each be in a different
infinite component.

Using the finite energy property to change the configuration $b$ inside of $B$
when it encounters three infinite thick components,
as illustrated on Figure~\ref{fi:Trifurcation},
there is a positive probability of observing a trifurcation
inside of $B$ for $\gamma_n(b)$.

The rest of the proof follows as in the original theorem:
if the probability that a cell is a trifurcation is positive,
then so is the frequency of trifurcations by Birkhoff's ergodic theorem on $\Z^d$,
thus it must be of order $n^d$ in a big hypercube.
However, a theoretical $O\left(n^{d-1}\right)$ bound can be obtained
on the amount of trifurcations, thus a contradiction.
The number $N$ cannot be infinite.

With a similar but much simpler finite energy argument, $N$ cannot be constant
greater or equal to $2$, as the probability of having at most $N-1$ components
would be positive, by opening an entire hypercube encountering several components.
\end{proof}
\end{lemma}

Thanks to this result,
we can from now on talk about \emph{the} infinite component of the percolation,
whenever it exists.
We now need to actually control the frequency of cells belonging to it.
Further analyses will be done on a Bernoulli noise,
but we still hope for a more general result to come from percolation theory.

\begin{proposition}[Frequency of the Infinite Component] \label{prop:DensityPerco}
Consider $I(b)\subset \Z^d$ the random infinite component of 
the $n$-thickened percolation $\gamma_n(b)$,
with respect to the original $\epsilon$-Bernoulli percolation 
$\mathds{P}=\B(\epsilon)^{\otimes\Z^d}$.

Then the constant $C_n^d= 48(2n+1)^d$ is such that
$\mathds{P}(0\notin I) \leq C_n^d\times\epsilon$.

\begin{proof}
Let us describe first what the event $\left\{0\notin I\right\}$ represents.
Either the cell $0$ is closed in $\gamma_n(b)$
(\emph{i.e.}\ $\gamma_n(b)_0=1$) so that it belongs to no component,
or it is open, but its component is finite.
The first scenario happens with probability $\left(1-(1-\epsilon)^{(2n+1)^d}\right)$.

In the second scenario, this implies that the component of $0$
in the sub-percolation induced by $\gamma_n(x)$
on the network $\Z^2\times \{0\}^{d-2}$ is also finite.
Consider the sub-network $[(2n+1) \Z]^2\times\{0\}^{d-2}$,
where two cells are adjacent whenever one coordinate differs by $2n+1$.
If two neighbouring cells $e$ and $f$ of this sub-network are open in $\gamma_n(b)$,
then all the cells in $\left(e+B_n\right)\cup\left(f+B_n\right)$ must be open.
Hence, if $e$ and $f$ are open, connected in the sub-network,
then all the cells that lie in-between in $\Z^2$ are also open,
so that $e$ and $f$ are in the same connected component of $\gamma_n(b)$.
The interest of this trick is that, as those windows $e+B_n$ and $f+B_n$ are
disjoint, the value of the cells $e$ and $f$ in $\gamma_n(b)$ are actually independent.
To put it short, in this second scenario,
the component of $0$ in the sub-network $[(2n+1)\Z]^2$ must be finite too.

The percolation on this sub-network is just a plane
$\left(1-(1-\epsilon)^{(2n+1)^d}\right)$-Bernoulli independent site percolation.
In this case, if the component of $0$ is finite,
then the outer boundary of this component must be a cycle of closed cells,
where two neighbouring cells may be diagonally adjacent,
so we just need an upper bound on the probability of \emph{this} event.

We can easily start with the upper bound $1-(1-\epsilon)^{(2n+1)^d}\leq 
(2n+1)^d\epsilon$ on the probability of a cell being closed.
Now we need to count the amount of cycles of a given length $l$.
Such a cycle must necessarily intersect the half-line $\N^*\times\{0\}$,
let's say at coordinates $(k,0)$,
and each of the columns $\{j\}\times \Z$ with $0\leq j <k$ must cross
the cycle at least twice, thus $l\geq 2k$ gives us an upper bound on the coordinate $k$.
Note also that a cycle is in particular a self-avoiding path,
so that, for a fixed value of $k$, we can upper bound the number of cycles by
$9\times 8^{l-1}$.
Whenever $\epsilon< \frac{1}{8(2n+1)^d}$, we have:
\[
\begin{array}{rcl}
\mathds{P}(0\notin I) &\leq & (2n+1)^d\epsilon + \sum\limits_{l\geq 4}
\frac{l}{2}\times 9\times 8^{l-1} \times \left((2n+1)^d\epsilon\right)^l \\
&\leq & \frac{9}{16}\epsilon \times \sum\limits_{l\geq 1}
8 (2n+1)^d\times l\left( 8(2n+1)^d \epsilon \right)^{l-1} \\
&= & \frac{9}{16}\epsilon \times\partial_\epsilon\left[
\sum\limits_{l\geq 0}\left( 8(2n+1)^d \epsilon \right)^l \right] \\
&=&\frac{9}{16}\epsilon\times
\partial_\epsilon\left[\frac{1}{1-8(2n+1)^d \epsilon}\right]
= \frac{9}{16}\epsilon\times \frac{8(2n+1)^d}{(1-8(2n+1)^d\epsilon)^2} . \\
\end{array}
\]

So far, this upper-bound is of the form $\epsilon f(\epsilon)$
for some function $f$ that is positive increasing on the interval
$\left[ 0,\frac{1}{8(2n+1)^d}\right[$ and goes to infinity on the right.
If we find $\epsilon_0$ in this interval such that
$\epsilon_0 f\left(\epsilon_0\right)=1$, then
the upper bound by $f\left(\epsilon_0\right) \epsilon$
will hold on this interval as $f$ is increasing,
and the upper bound will hold for $\epsilon_0 \leq \epsilon \leq 1$ as
$\mathds{P}(0\notin I)\leq 1 \leq f\left(\epsilon_0\right) \epsilon$ on this interval.

Let us denote $a=\frac{9}{16}$ and $b=8(2n+1)^d$.
Solving $\epsilon f(\epsilon)=1$ equates finding the root
of $b^2 \epsilon^2 -b(a+2)\epsilon+1$ on the interval $\left[ 0,\frac{1}{b}\right]$.
The roots are $\epsilon_\pm=\frac{\sqrt{a+2}}{b}
\left(\frac{\sqrt{a}+\sqrt{a+2}}{2}\right)$
and only $ \epsilon_-$ is in the desired interval.
A direct computation then yields $f\left(\epsilon_-\right)=
\frac{2b}{1-a\left(\sqrt{1+\frac{2}{a}}-1\right)}$.
Replacing $a$ by its value, we obtain
$1-a\left(\sqrt{1+\frac{2}{a}}-1\right)=\frac{25-3\sqrt{41}}{16}>\frac{1}{3}$,
thus finally $f\left(\epsilon_-\right)< 6b$.
At last, the constant $C_n^d= 48(2n+1)^d$ provides the desired upper bound.
 \end{proof}
\end{proposition}

This proof depends on the specific properties of the independent percolation,
but is quite elementary in exchange.
In order to adapt the following periodic stability theorem
to a more general class of $\epsilon$-noises,
one would first need to obtain a similar lower bound
on the frequency of cells in the infinite connected component,
the equicontinuity of $\mathds{P}(0\notin I)$ as $\epsilon\to 0$.

\subsection{Periodic Stability Theorem}

\begin{theorem}[Periodic Stability Theorem] \label{thm:PeriodicStability}

Consider $\Omega_\F$ a 2D+ periodic SFT.
Then $\Omega_\F$ is $f$-stable for $d_B$ on the class $\B$ of Bernoulli noises,
with linear speed $f(\epsilon)=2C_{c(\F)}^d\epsilon$.

\begin{proof}
In order to obtain linear stability, we will consider a measure
$\lambda\in\widetilde{\M_\F^\B}(\epsilon)$, and build
a measurable mapping $\psi:\Omega_{\widetilde{\F}}\to\Omega_\F$,
so that $d_H\left(\omega,\psi(\omega,b)\right)$ is small
for a $\lambda$-typical configuration $(\omega,b)\in \Omega_{\widetilde{\F}}$.

Consider $N$ the size of a base hypercube for the periodic SFT $\Omega_\F$,
and $c$ the constant obtained in Lemma~\ref{lem:LocalGlobal}.
As $\A^{\llbracket 0, N-1\rrbracket^d}$ is finite, then so is $\Omega_\F$.
Thus, it makes sense to consider $\Omega_\F$ as a finite alphabet
and to define the full-shift $\Omega_{\Omega_\F}$.

Let us define the morphism
$\rho:\Omega_{\widetilde{\F}}\to \Omega_{\Omega_\F}$
such that, whenever the window $B_c$ is clear in
$\sigma_e(\omega,b)\in\Omega_{\widetilde{\F}}$,
then $\rho(\omega,b)_e=\omega_e$ as in Lemma~\ref{lem:LocalGlobal},
but specifically for the window $B_c$ of $\sigma_e(\omega,b)$ centred on $0$.
If the window is obscured, then we may default to
some configuration $\omega'\in\Omega_\F$.
The interest of ``forgetting'' the role of the coordinate $e$,
of acting as if each cell was the centre of the network $0\in\Z^d$,
is that this way $\rho$ is $\sigma$-invariant,
we have a local characterisation of the morphism
$\rho : \widetilde{\A}^{B_c} \to \Omega_\F$.

Without loss of generality,
assume the finite set $\left(\Omega_\F,<\right)$ is strictly ordered.
We may now define the adjusted majority rule cellular automaton
$\theta_n: \Omega_\F^{B_n} \to \Omega_\F$ as follows.
First, map each configuration of the pattern $\left(\omega_e\right)_{e\in B_n}$ 
onto the configuration $\sigma_{-e}\left(\omega_e\right)$,
so that we locally undo the offset introduced by $\rho$
by aligning all the configurations on a ``common'' centre.
Only then we may apply a regular majority rule,
on the family $\left(\sigma_{-e}\left(\omega_e\right)\right)_{e\in B_n}$,
by picking the maximal configuration for the arbitrarily introduced order in case of a tie.

Consider now the morphisms $\psi_n=\theta_n\circ\rho$
obtained by applying an adjusted majority rule over $\rho$.
Using once again the order on $\Omega_\F$, we may define the pointwise limit
$\psi=\varlimsup\limits_{n\to\infty} \psi_n$,
which is still $\sigma$-invariant and measurable.
Note that the value of $\psi(\omega,b)$ in some cell may now depend on arbitrarily
far values, so $\psi$ is \emph{not} a morphism.

Consider the configuration $(\omega,b)\in \Omega_{\widetilde{\F}}$, and let
$I\subset \Z^d$ be the infinite component of the $c$-thickened percolation in $b$.
As $\omega|_{I+B_c}$ is locally admissible,
$\omega|_I$ is a globally admissible pattern,
the restriction of some configuration $\omega_0\in\Omega_\F$.
For any cell $e\in I$, we have $\rho(\omega,b)_e=\sigma_e\left(\omega_0\right)$.

Assume now that $\epsilon<\frac{1}{2C_c^d}$,
so that in the Bernoulli percolation, $I$ has a density greater than $\frac{1}{2}$
according to Proposition~\ref{prop:DensityPerco}.
This means that, $\lambda$-almost-surely,
after some rank $n_0$, strictly more than half of the cells $f\in e+B_n$ 
of $(\omega,b)$ are inside of $I$,
thus are mapped by $\rho$ onto translations $\sigma_f\left(\omega_0\right)$.
Thence, after the very same rank $n_0$,
$\psi_n(\omega,b)_e=\sigma_e\left(\omega_0\right)$.
Consequently, by taking the limit $n\to\infty$,
$\lambda$-almost-surely,
$\psi(\omega,b)_e= \sigma_e\left(\omega_0\right)$ for any cell $e\in\Z^d$.

To sum it up,
$(\omega,b)\mapsto \psi(\omega,b)_0=\omega_0$
is a measurable mapping $\Omega_{\widetilde{\F}}\to\Omega_\F$,
such that $d_H\left( \omega,\omega_0\right)\leq C_c^d \epsilon$
whenever $\epsilon\leq \frac{1}{2 C_c^d}$.
More generally,
the bound $d_H\left(\omega,\psi(\omega,b)_0\right)\leq 2 C_c^d \epsilon$
holds $\lambda$-almost-surely for any choice of $\epsilon$,
which finally gives us the linear bound we wanted:
\[
d_B\left(\pi_1^*(\lambda),\M_\F\right) \leq
d_B\left(\pi_1^*(\lambda),\left[\psi(\cdot)_0\right]^*(\lambda)\right)
\leq 2 C_c^d \epsilon .
\]
\end{proof}
\end{theorem}

This concludes our analysis of periodic SFTs in the 2D+ case.
The explicit constant $C_n^d$ could doubtlessly be improved,
but such matters would require much more work without improving on the
\emph{linear} aspect of the bound.

A further track of reflection, as already mentioned earlier, may be to extend
this theorem to a more general class of noises, using stronger percolation results,
while leaving much of the actual proof of the theorem unchanged.

What we got interested in instead is the study of stability for aperiodic SFTs.
We chose the well-known Robinson tiling, as it is already \emph{almost} periodic,
in order to adapt the previous scheme of proof as much as possible.
This will be the topic of the last section of the paper.

\section[The Case of 2D (c1,c2)-Robust Tilings]{The Case of 2D $(c_1,c_2)$-Robust Tilings} \label{sec:RobustDurRoShe}

Before diving into the Robinson tiling, let us now digress a bit
to contextualise our study.
The aim of this section is to provide an informal analysis
of an already existing Besicovitch stability result in our current framework.
More precisely, we are interested in the notion of stability
described by Durand, Romaschenko and Shen~\cite{DuRoShe12},
which was then used to prove periodic stability in the 2D case
in a further article by Ballier, Durand and Jeandel~\cite{BaDuJean10}.

Here, we will provide a rough and qualitative estimate of
the convergence speed obtained with their method.
Yet, for this article to be as self-contained as possible,
we will still introduce the essential definitions to understand the cited results.

The estimates provided here bear no influence on the following aperiodic stability result,
so this section can easily be skipped in a first reading of the current article.

\subsection{Robust Tilings and Sparse Sets}

To obtain stability, instead of using a notion of percolation
-- which is best seen as a clear connected tree that spans the whole obscured space --
they introduce the notion of islands of errors
-- which is best seen as small clumps of obscured cells isolated in the whole clear space.

\begin{definition}[$(\alpha,\beta)$-Island of Errors]
Consider a noise configuration $b \in \{0,1\}^{\Z^2}$
which we identify to $E\subset \Z^2$ the set of obscured cells.

A set $F\subset E$ is an $(\alpha,\beta)$-island of $E$ if
$F$ can be included in some $\alpha$-square
and its $\beta$-neighbourhood does not meet any other obscured cell of $E$,
\emph{i.e.}\ $\left(F+B_\beta \right) \cap \left( E\backslash F\right)=\emptyset$.
\end{definition}

In this framework, the ``right'' way to obtain stability
is to remove the islands of obscured cells,
by changing the values of the tiles underneath on a small neighbourhood.
This is well-encapsulated by the following notion of robustness.

\begin{definition}[$\left(c_1,c_2\right)$-Robustness]
Let us denote by $R_{i,j}=S_j\backslash S_i$ (with $i<j$) the ring-shaped window
obtained by removing the $i$-square $S_i$ at the centre of a $j$-square.

Let $0<c_1\leq c_2$ be two positive integers.
A set $\F$ is $(c_1,c_2)$-robust if,
for any $n\in\N$ and any locally admissible pattern $u\in \A^{R_{n,c_2 n}}$,
there exists a locally admissible pattern $v\in \A^{S_{c_2 n}}$
such that $u$ and $v$ coincide on $R_{c_1 n,c_2 n}$
-- which is a strict subset of the ring $R_{n,c_2 n}$ as long as $c_1\geq 2$.
\end{definition}

An explicit example of robust tiling is any periodic SFT~\cite{BaDuJean10},
roughly for the same reason we could obtain a globally admissible configuration
by peeling a constant width of the border of any pattern in the previous section.
However, this notion is much more general,
and strongly aperiodic robust SFTs are proven to exist~\cite{DuRoShe12}.

Note that, while the constants may change in the process,
this notion of robustness is stable under conjugacy,
so that we cannot prove stability of a non-robust SFT by looking
for a suitable robust conjugated.

Whenever $\beta \geq c_2 \alpha$, we are roughly in a situation where we can ``repair''
an island of errors by changing the tiles in a $c_1 \alpha$-square.
Hence, we need some guarantees that $E$ is entirely made out of islands we can correct.

\begin{definition}[$(\alpha,\beta)$-Sparse Set]
A set $E=E_0$ is said to be sparse,
given a sequence $\left(\alpha_k,\beta_k\right)_{k\in\N^*}$,
if we can step by step remove all the $\left(\alpha_k,\beta_k\right)$-islands
from $E_{k-1}$ to obtain a set $E_k$,
in such a way that the decreasing limit set $E_\infty=\bigcap E_k$ is empty.
\end{definition}

Up to now, the definitions introduced were formal.
For the rest of this section, we will provide
a qualitative and quite handwavy analysis
of the convergence speed we can obtain in this framework.

\subsection{Qualitative Convergence Speed}

By the Borel-Cantelli theorem, any $\epsilon$-Bernoulli noise
will certainly contain islands for any pair $(\alpha,\beta)$,
which may \emph{a priori} be hard to correct.
However, it is proven~\cite[Lemma 3]{DuRoShe12} that,
assuming $8 \sum_{k=1}^{n-1} \beta_k < \alpha_n \leq \beta_n$ for any $n\in\N^*$
and $\sum_n \frac{\ln\left(\beta_n\right)}{2^n}<\infty$,
then for $\epsilon$ small enough the random set $E$ is
almost-certainly $(\alpha,\beta)$-sparse.
Unfortunately, general bounds on $\epsilon$ would be quite hard to obtain,
but we will provide rough estimates for our choice of $(\alpha,\beta)$.

It is also proven~\cite[Lemma 4]{DuRoShe12} that
in any $(\alpha,\beta)$-sparse set $E$,
the density of obscured cells is \emph{at most}
$\sum_n \left(\alpha_n / \beta_n\right)^2$
-- the main argument is that each $\left(\alpha_k,\beta_k\right)$-island
contains \emph{at most} $\alpha_k^2$ obscured cells,
among \emph{at least} $\beta_k^2$ cells in a neighbourhood of the island
disjoint of the other islands and their neighbourhoods.
To properly quantify the convergence speed,
we would need to take into account the density not
of the islands of errors but of the $c_1\alpha$-square around them,
but this estimate will suffice for the present qualitative analysis.

Note that, as $\frac{\alpha_n}{\beta_n}$ must go to $0$ for this sum to be finite,
as we will have to take $k\to \infty$ as $\epsilon\to 0$
for the density of errors to vanish, then naturally
a bound on the convergence speed obtained by this method will hold true
for \emph{any} pair $\left(c_1,c_2\right)$.
The tricky part will be that the domain $[0,\tau]$ itself
on which the bound holds will depend on the pair.

Consider $\alpha_n = 8^n (n-1)!n!$ and $\beta_n=8^n (n!)^2$.
It is clear that any $k$-\emph{shift} of this sequence
(starting at some rank $k+1$ instead of $1$)
will satisfy the previously stated hypotheses.
For a given sparse set $E$ for the $k$-shifted sequence,
the density of errors is $\sum_{n=k+1}^\infty\frac{1}{n^2}
\leq \int_k^\infty \frac{1}{t^2}\mathrm{d}t = \frac{1}{k}$.

To obtain the convergence speed,
we now need to estimate the maximal value of $k$ such that $E$ is sparse
for the $k$-shifted sequence for a given $\epsilon$.
Looking at the proof of the result~\cite[Lemma 3]{DuRoShe12},
it appears that the key property to obtain sparsity is that
$\sum_n \frac{\ln\left(\beta_n\right)}{2^n} < \ln\left(\frac{1}{\epsilon}\right)$.
As $\ln\left(\beta_n\right)= 8\ln(n)+2\ln(n!) \leq n^2$ after some rank,
for the $k$-shifted sequence, we can bound the left term by $k^2+4k+6$.
Asymptotically, the best choice for $k$ is thus $k \approx \sqrt{\ln(1/\epsilon)}$,
so that $f(\epsilon)\approx \frac{1}{\sqrt{\ln(1/\epsilon)}}$.

Considering all the small approximations we did on the way,
what matters here is not the value of the bound but its order of magnitude.
Indeed, $\frac{1}{\sqrt{\ln(1/\epsilon)}}$
is much \emph{much} slower than any polynomial speed,
which legitimises our efforts to obtain
a linear convergence speed in the periodic case.

The notion of islands and sparsity can be used as a black box
to obtain percolation results~\cite[Section 9.3]{DuRoShe12},
hence as a tool it is in some ways more powerful
than the percolation theory we used in the previous section.
However, as we have seen here, this versatility comes at the cost
of the precision and simplicity of the bounds we can obtain.

\section{The Robinson Tiling: an Almost Periodic Stable Example} \label{sec:Robinson}

The first aperiodic tiling defined by local rules was
proposed by R. Berger~\cite{Ber66},
who used 20426 Wang tiles to encode a hierarchical structure,
and thus aperiodicity.
The construction was strongly simplified by R. Robinson~\cite{Rob71}
who proposed a Wang tileset with 56 tiles,
which once again forces a hierarchical structure.
In fact, if we allow diagonal interactions between tiles,
the number of tiles can be brought down to 6 tiles
and their rotations and symmetries~\cite{Rob71}.
The simplicity of the tileset and its hierarchical structure,
with arbitrary large squares which permits the embedding of
space-time diagrams of Turing machines into it,
explains why the Robinson tiling is certainly the most studied aperiodic tiling.

The Robinson tiling is \emph{not} $(c_1,c_2)$-robust in the previous sense:
it can have an infinite central cross in $\Z^2$ with a black arm in each direction,
with only one obscured cell at the centre,
that no amount of removing may correct.
Thus, it seems difficult to correct mistakes locally here.
However, the hierarchical structure implies that for a given scale,
the corresponding squares form a periodic structure,
except for a small fraction of tiles that corresponds to
the squares higher in the hierarchy.
A similar technique that in Section~\ref{sec:2DPeriodic}
yields some stability at this scale,
and allows us to deduce the stability of the Robinson tiling
with a polynomial speed (Therorem~\ref{thm:RobinsonStability}).

\subsection{The Classic Robinson Tiling}

\begin{figure}[H]
\centering
\includegraphics[scale=.5]{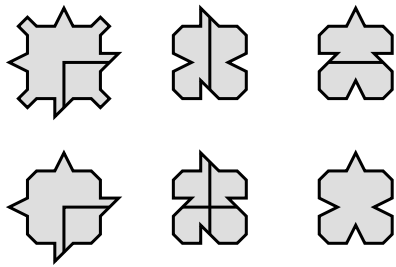}
\caption{The six base Robinson tiles.}
\label{fi:VanillaRobinson}
\end{figure}

Our first attempt at 2D+ aperiodic stability used
the folkloric Robinson tiles shown in Figure~\ref{fi:VanillaRobinson},
and their rotations and symmetries --
so that the total amount of tiles is actually $32$.

With this tileset, the forbidden patterns are self-evident:
two laterally adjacent tiles must have matching borders,
including the black lines drawn on them,
and any square made of four tiles must use exactly one rotation
of the top-left tile in Figure~\ref{fi:VanillaRobinson} with bumpy corners,
so that the small diamond in the centre of the square is filled-in.
Any non-matching pair or square of adjacent tiles is then a forbidden pattern in $\F$.

Note that forbidden patterns can occur with left-right and top-bottom neighbours,
but \emph{also} on diagonally adjacent tiles,
unlike the tiling originally introduced by Robinson in the context of Wang tiles.
The two tilings can still be easily conjugated.

\begin{figure}[H]
\centering
\includegraphics[scale=.35]{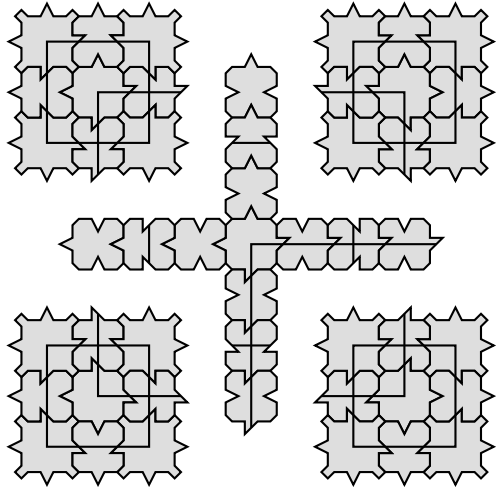}
\caption{Four $2$-macro-tiles around a central cross form a $3$-macro-tile.}
\label{fi:MacroTile}
\end{figure}

\begin{definition}[Macro-tiles]
We define macro-tiles inductively.
First, the $1$-macro-tile is just the top-left tile of Figure~\ref{fi:VanillaRobinson},
with bumpy corners.

Then, the $(N+1)$-macro-tile is obtained by sticking four $N$-macro-tiles
in order to draw a square around a central cross,
as shown in Figure~\ref{fi:MacroTile}.
\end{definition}

\begin{definition}[Orientation Symbols]
Let us use the symbol \rbtile{} to denote
the \emph{default} orientation of a $N$-macro-tile,
with the black arms of the central cross pointing on the bottom and on the right,
as seen in Figure~\ref{fi:MacroTile}.
Likewise we denote \lbtile{}, \lttile{} and \rttile{} for the other orientations.
\end{definition}

By induction, we have that a $N$-macro-tile is a $2^N-1$ tiles long square.
One can prove that two $N$-macro-tiles cannot overlap.
These fundamental properties can be found in Robinson's seminal article
\emph{Undecidability and nonperiodicity for tilings of the plane}~\cite{Rob71},
and are nicely condensed into seminar notes~\cite{Schwa07}.

A Robinson tiling is almost periodic,
in the sense that any given window of a tiling occurs periodically in the tiling,
but not always with the same periodicity.
Most notably, if you keep only the $N$-macro-tiles and forget about
the thin grid sticking all of them together,
you obtain a $2^{N+1}$-periodic pattern,
which has density $\left(1-\frac{1}{2^N}\right)^2$.

\begin{figure}[H]
\centering
\includegraphics[scale=.4]{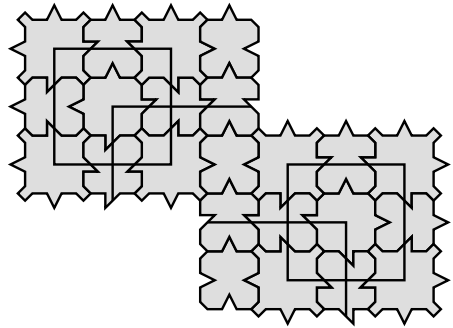}
\caption{Two loosely aligned $2$-macro-tiles,
with one tile in common and a tiled gap.}
\label{fi:Align}
\end{figure}

The issue with this tiling is that the alignment of macro-tiles on such a grid
is a consequence of the global structure of a Robinson tiling,
and is not enforced by the local rules.
This is illustrated by the two misaligned macro-tiles in Figure~\ref{fi:Align},
and such a phenomenon can arise at \emph{any} scale.
This implies that we would not be able to ensure stability using a percolation argument
as we did for the periodic case.

By pushing this phenomenon to the limit,
we can obtain ``pathological'' Robinson tilings that exhibit a \emph{cut},
an infinite horizontal or vertical line,
with a misalignment on both sides.

\subsection{An Enhanced Robinson Tiling}

To work around the aforementioned issue,
let us now introduce a variant tileset
by \emph{adding} information over the already existing tiles.

To force this alignment in a local way,
we want for each macro-tile to send a ``signal'' from its central cross,
which will force the correct alignment between neighbouring macro-tiles at any scale.
The idea originates in Sylvère Gangloff's phd thesis~\cite{Gan18}
on another variant of the Robinson tiles,
and we transpose it on our current tileset.

\begin{figure}[H]
\centering
\includegraphics[scale=.5]{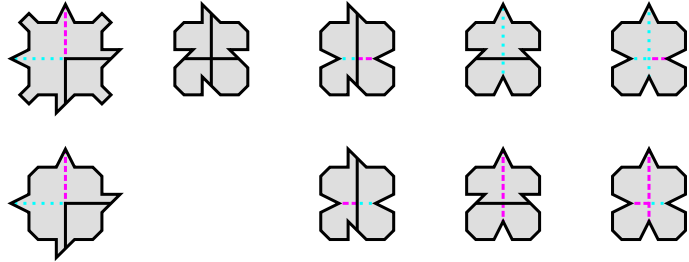} 
\caption{The nine enhanced Robinson tiles.}
\label{fi:SelfAlignRobinson}
\end{figure}

More precisely, consider the tiles on Figure~\ref{fi:SelfAlignRobinson},
roughly grouped according to which of the previous tiles they come from.
Now, all of the tile have a cross-like pattern drawn upon them.
In order to preserve their specific orientation,
the two leftmost tiles must never undergo a symmetry,
so that a \rbtile{} tile always has a blue dotted line pointing left
and a red dashed one pointing up.
Up to symmetry of the other tiles and rotation,
this brings the total amount of tiles to 56.

We define the set of forbidden patterns as before,
now in accordance with the crosses drawn upon the tiles.
For the rest of the section, we will use $\F$
to denote this specific set of forbidden patterns.
Using the same process as before, starting from the base $1$-macro-tiles,
there is a unique way to build macro-tiles inductively.
For a given macro-tile of the initial Robinson tiling,
we can without ambiguity deduce where the red dashed lines and blue dotted lines
of the enhanced macro-tile are.

As there is a direct local projection (thus a morphism)
of this enhanced tileset on the previous Robinson tiles,
any configuration is still aperiodic.
However, this morphism is not a bijection.
On one hand, this morphism is not surjective, as we cannot reach
tilings with a misaligned cut.
On the other hand, this morphism is not injective,
as we may have an aligned cut with an infinite red dashed \emph{or} blue dotted line
that gets projected onto the same configuration.
The main interest of this added structure, as we will prove, is that
it indeed \emph{locally} enforces the alignment we lacked before.

\subsection{Local Alignment Properties}

As we already said, we want to study the almost periodicity
obtained by looking only at $N$-macro-tiles.

\begin{definition}[Well-Aligned and Well-Oriented Pairs]
A pair of $N$-macro-tiles (here seen as a pattern in $\Z^2$) is said to be 
well-aligned if both of their centres have one coordinate in common,
and the other differs by exactly $2^N$ so that there is a gap
of precisely one line/column between them.

More generally, we say the two $N$-macro-tiles are
\emph{loosely} aligned (with $0<k\leq 2^N-1$ tiles in common)
when one of the coordinates of their centres differs by exactly $2^N$
and the other by $2^N-k-1$,
\emph{i.e.}\ we start with a well-aligned pair (with $2^N-1$ tiles in common)
and we translate one of them of $k$ units in the direction of the gap in-between.

A pair of well-aligned macro-tiles is said to be well-oriented if their central crosses
form a pattern \rbtile{} \lbtile{} or \lbtile{} \rbtile{} (or a rotation of these),
which can actually be filled by a central cross
in the process of making a larger macro-tile.
\end{definition}

\begin{definition}[Edge Words of Macro-Tiles]
We define the words $l_N$ and $t_N$,
obtained by reading the colours on the left and top edges
of the \rbtile{} $N$-macro-tile in a clockwise motion,
with blue dotted lines encoded as a $0$ and red dashed lines as a $1$.
\end{definition}

For a binary word,
we define $\overline{b}=1-b$ the binary complement of a letter,
extended to binary words by a direct induction.
We also define the $\text{mirror}$ function on words
such that $\text{mirror}(uv)=\text{mirror}(v)\text{mirror}(u)$,
that returns the same word but backwards.
Both of these mappings are involutions
and they commute with each other.

\begin{lemma}
For any $N\in\N^*$,
we have $t_N=\overline{\text{mirror}\left(l_N\right)}$.

What is more, $\left|l_N\right|=\left|t_N\right|=2^N-1$ is odd,
and these words actually differ of exactly one letter in their middle.

\begin{proof}
For $N=1$, we simply have $l_1=0$ and $t_1=1$.

By induction, as seen in Figure~\ref{fi:MacroTile},
when building a $(N+1)$-macro-tile, on the left half from bottom to top,
we first have a \rttile{} $N$-macro-tile that reads as $t_N$,
then we read the $0$ given by the blue dotted arm of the central cross,
and finally $l_N$ on the \rbtile{}, so that $l_{N+1}=t_N 0 l_N$.
Likewise, $t_{N+1}=t_N 1 l_N$.
Hence:
\[
\overline{\text{mirror}\left(l_{N+1}\right)}=
\overline{\text{mirror}\left(t_N 0 l_N\right)}=
\overline{\text{mirror}\left(l_N\right)}1
\overline{\text{mirror}\left(t_N\right)}=
t_N1l_N = t_{N+1} ,
\]
which concludes the proof by induction.
\end{proof}
\end{lemma}

\begin{figure}[H]
\centering
\includegraphics[scale=.5]{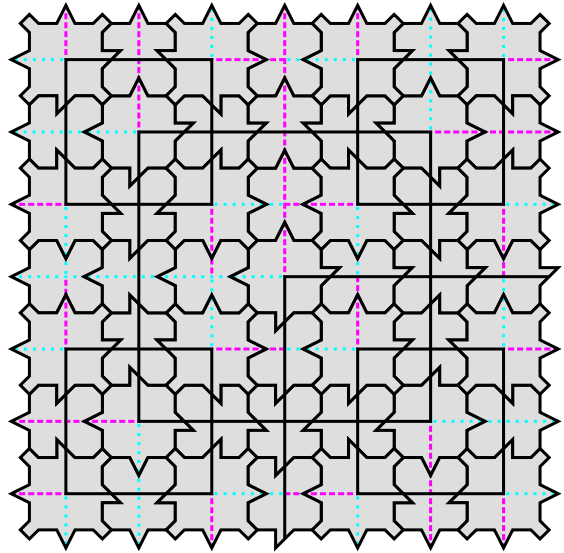} 
\caption{The $3$-macro-tile obtained using the enhanced tileset.}
\label{fi:EnhancedMacro}
\end{figure}

In Figure~\ref{fi:EnhancedMacro}, for example,
we observe that $l_3=1100100$ and $t_3=1101100$.

\begin{proposition}[Local Alignment of Macro-Tiles] \label{prop:LooselyAligned}

For any scale $N\in\N^*$, a pair of loosely aligned $N$-macro-tiles
with a tileable gap in-between \emph{must} be well-aligned and well-oriented.

\begin{proof}
Assume first that two well-aligned macro-tiles are not well-oriented.
If only one of these tiles has a black arm that falls into the gap
(\emph{e.g.}\ a \rbtile{} \rttile{} pattern), then this gap cannot be tiled.
Up to a rotation, the remaining cases are the
\lbtile{} \rttile{} and \lttile{} \rbtile{} patterns.
In these cases, the right arm of the left cross and the left arm of the right cross
have the same colour, thus no tile can fill the gap in-between.
In other words, by contraposition,
a well-aligned pair with a tileable gap must be well-oriented.

At the scale $1$, if two tiles are loosely aligned they are actually well-aligned,
thus if the gap is tileable they are well-oriented.
This allows us to initialise the induction.

Assume the result holds up to scale $N\in\N^*$
and consider a pair of $(N+1)$-macro-tiles,
once again loosely aligned with a tileable gap.
The macro-tiles cannot have exactly one tile in common,
which would imply that we have two $1$-macro-tiles well-aligned
with a tileable gap but ill-oriented, hence $k\geq 2$.

What is more, $k$ cannot be even.
Assuming $k$ is even, this pair of $(N+1)$-macro-tiles
contains a pair of $2$-macro-tiles with a tileable gap and $2$ tiles in common.
It is clear that this cannot happen, by an exhaustion of cases.
For example, looking at a well-aligned \lbtile{}~\rbtile{} pair,
if we move the right tile of one unit upwards,
then the right arm of the left tile and the bottom-left corner of the right tile
face a tileable gap with a red dashed line, which is impossible.

This concludes the case $N+1=2$,
as $k\geq 3$ must then be equal to $3$, maximal,
so that the $2$-macro-tiles are well-aligned.
Likewise, when $N+1>2$,
the $N$-macro-tiles must be well-aligned with $k$ odd,
so either the $(N+1)$-tiles are well-aligned,
or only half of the $N$-macro-tiles actually face the gap and are well-aligned.
In the second scenario, we are once again in a tileable ill-oriented case, impossible.
Finally, the $(N+1)$-macro-tiles must be well-aligned thus well-oriented,
which concludes the induction.
\end{proof}
\end{proposition}

\begin{proposition}
For any scale $N\geq 2$, consider the constant $C_N=2^N-1$,
such that for any $n\in\N$ and
any clear locally admissible pattern $\omega$ on $B_{n+C_N}$,
its restriction $\omega|_{B_n}$ is almost globally admissible,
in the sense that up to a low-density grid,
$\omega|_{B_n}$ is the restriction of a Robinson tiling,
with well-aligned and well-oriented $N$-macro-tiles.

\begin{proof}
We will demonstrate a slightly stronger result here, \emph{i.e.}\
that by removing \emph{at most} $C_N$ layers of tiles
on the top, bottom, left and right sides of any locally admissible square,
and not necessarily the same amount of layers on each side,
we obtain an actual family of well-aligned and well-oriented $N$-macro-tiles
with respect to their neighbours.
Thence, by actually peeling $C_N$ layers on each side,
we obtain the stated result.

To do so, we need to proceed inductively, as before.
We cannot initialise the result at $N=1$,
but notice that if the result holds at rank $N$ with the constant $C_N$,
then it also holds at any lower rank with the same constant.

\begin{figure}[H]
\centering
\includegraphics[width=\linewidth]{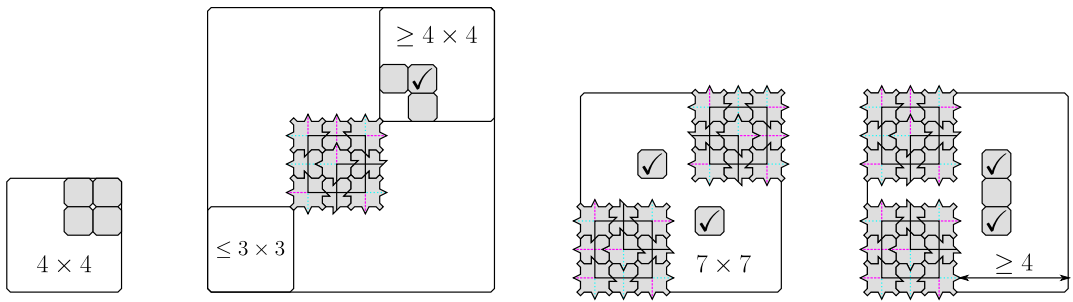}
\caption{From left to right, key steps A, B, C and D of the case $N=2$.}
\label{fi:Proof}
\end{figure}

Hence, let us now prove the case $N=2$ for a tiled $n$-square $B$.
First, it is known that with the initial Robinson tileset,
if a $3$-square is tiled with a \rttile{} in the bottom-left corner,
then it is tiled by a $2$-macro-tile.
This property still holds for the enhanced tileset,
and can be easily checked by enumerating all the cases.

We will inductively build a rectangle of well-aligned $2$-macro-tiles
in the $n$-square $B$, assuming that $n\geq 10$ for now.
If we look at the $4$-square in the bottom-left corner,
one of the four cells highlighted in the step A of Figure~\ref{fi:Proof}
must contain a $1$-macro-tile, with bumpy corners.
Then this bumpy corner is actually part of a $2$-macro-tile in $B$.
So far, we have a $1\times 1$ rectangle of $2$-macro-tiles.

As illustrated in step B,
considering where our first bumpy corner was,
there is at \emph{most} a $3\times 3$ rectangle in the bottom-left corner
(diagonally adjacent to the $2$-macro-tile), and $n\geq 10$,
so the top-right corner is at \emph{least} a $4\times 4$ rectangle of tiles.
One of the three highlighted tiles in step B must be a bumpy corner too.
If there was a corner in one of the unchecked tiles,
it would be part of a $2$-macro-tile,
that should either intersect the one drawn on Figure~\ref{fi:Proof}
-- which is impossible even for regular Robinson tiles -- or be
loosely aligned with it -- which is impossible
according to Proposition~\ref{prop:LooselyAligned}.
Hence the checked cell must contain a tile with bumpy corners,
and more precisely a $\rttile{}$ for the same reasons.
This tile can then be completed into a $2$-macro-tile, which brings us to step C.
There, the two checked cells must contain a $1$-macro-tile too,
and each can be completed into its own $2$-macro-tile,
so that we obtain at last a $2\times 2$ rectangle of $2$-macro-tiles.

Just like the two diagonally adjacent $2$-macro-tiles present in step C imply
a square of $2$-macro-tiles,
the presence of two laterally adjacent $2$-macro-tiles in step D
implies a square of $2$-macro-tiles.
Thus, now that we have a rectangle with \emph{at least} $2$ macro-tiles
on each side, we can repeat step D in each direction
as long as $4$ tiles or more remain.
Hence, as long as $n\geq 10$, $C_2=3$ works well.

More largely, we can entirely peel a $9$-square if we remove $5$ layers
on each side, so we proved that $C_2=5$ works.
However, a more careful study of the cases $ n\in\{7,8,9\}$
allows us to conclude that $C_2=3$ works for these cases and is optimal
(to do so consider a $9$-square \emph{centred} on a $2$-macro-tile,
so that all the adjacent ones will be missing a layer).
When $n\leq 6$, $C_2=3$ trivially works too, which concludes our study of $N=2$.

Assume now that the result holds at rank $N$ with the constant $C_N$
and let us prove it at rank $N+1$.
We can start by peeling away at most $C_N$ tiles, using our induction hypothesis,
to obtain a grid of well-aligned and well-oriented $N$-macro-tiles.
A square of well-aligned $N$-macro-tiles can either form one $(N+1)$-macro-tile,
represent the lateral interface between two $(N+1)$-macro-tiles
or represent the central corner between four $(N+1)$-macro-tiles.
Thus, by peeling at most one layer of $N$-macro-tiles on each border --
a $N$-macro-tile not part of a $(N+1)$-macro-tile and the following grid,
so $2^N$ tiles in total -- we remove the incomplete interfaces and corners
to obtain a grid of well-aligned
$(N+1)$-macro-tiles (hence well-oriented by the previous proposition).
In conclusion, the result holds at rank $N+1$
with the constant $C_{N+1}=C_N+2^N$,
hence $C_N=2^N-1$ by a direct induction.
\end{proof}
\end{proposition}

\subsection{Almost-Stability at a Fixed Scale and Stability}

\begin{proposition}[Almost-Stability]
Let $\Omega_\F$ be the enhanced Robinson tiling.
For any choice of $\epsilon>0$, any scale $N\in\N^*$,
and any measure $\mu\in\M_\F^\B(\epsilon)$,
we have a coupling that yields:
\[
d_B\left(\mu,\M_\F\right)\leq 96\left(2^{N+2}+1\right)^2\epsilon
+ \frac{1}{2^{N-1}}.
\]

\begin{proof}
For a given scale $N$,
we want to apply the percolation argument
as if we were looking at a $\left(2\times 2^N\right)$-periodic SFT.
This added factor $2$ comes from the fact that,
for any globally admissible configuration,
the $(N+1)$-macro-tiles are well-aligned on a grid, and indistinguishable if we ignore
their central cross, hence the $N$-macro-tiles form a unique $2^{N+1}$-periodic pattern
up to translation.

By looking at a \emph{globally admissible} $2^N$-square,
we can always identify one, two or four partial $N$-macro-tiles
arranged in a square pattern around a central cross.
Thus, we can actually identify to which translation of the $2^{N+1}$-periodic pattern
this window corresponds.
Note that unlike in the general $k$-periodic case,
where we needed to look at $k$-squares to identify the translation,
we only need to look at a window of size $\frac{k}{2}$ here because
the Robinson tiling has a lot of intrinsic redundancy.

Just like in the periodic case, we can then look at the $c$-thickened percolation,
with $c=\left\lceil \frac{2^{N+1}+1}{2}\right\rceil + C_N=2^N+1+2^N-1=2^{N+1}$,
as explained in Lemma~\ref{lem:LocalGlobal}.
As stated in Proposition~\ref{prop:DensityPerco},
the infinite component of the $c$-thickened percolation
has density at least $1-48\left(2c+1\right)^2\epsilon$.

Let us add a blank symbol $\square\notin\A$ to the original alphabet.
Then, following the proof of Theorem~\ref{thm:PeriodicStability},
we can measurably map a noisy configuration $(\omega,b)$ onto
a globally admissible configuration $\psi(\omega,b)\in \Omega_\F$
but on the extended alphabet $\A\sqcup\{\square\}$,
such that almost-surely:
\[
d_H\left(\omega,\psi(\omega,b)\right) \leq 96\left(2^{N+2}+1\right)^2\epsilon
+ \frac{2^{N+1}-1}{2^{2N}} .
\]
The second term comes from the density of the symbols $\square$ in $\psi(\omega,b)$,
of the grid itself, which is equal to $1-\left(\frac{2^N-1}{2^N}\right)^2$.

In order to conclude, we need to explain how to measurably project $\psi(\omega,b)$
back onto the original alphabet $\A$, how to fill-in the grid,
so that we obtain an actual globally admissible enhanced Robinson tiling.
To do so, we can simply consider some measure $\widetilde{\mu}\in\M_\F$,
take a configuration $y\in\Omega_\F$ at random
independently of the rest following $\widetilde{\mu}$,
and then replace $\psi(\omega,b)$ by $\psi'(\omega,b,y)$ which is the unique translation
of $y$ by a vector $k\in\left\llbracket 0, 2^{N+1}-1\right\rrbracket^2$
such that the $N$-macro-tiles of $\psi'(\omega,b,y)$ and $\psi(\omega,b)$ are aligned.
This whole process is measurable, $\sigma$-invariant,
and only changes the values of $\psi(\omega,b)$ on the $\square$ tiles
which were already taken into account in the upper bound,
so that the same bound holds for $d_H\left(\omega,\psi'(\omega,b,y)\right)$.

Thence, we have a coupling such that
$d_B\left(\mu,\M_\F\right) \leq 96\left(2^{N+2}+1\right)^2\epsilon
+ \frac{1}{2^{N-1}}$, which proves the bound.
\end{proof}
\end{proposition}

By taking $N$ arbitrarily large, and then $\epsilon\to 0$,
we directly deduce the \emph{stability} of our enhanced Robinson tiling
for the Besicovitch distance.
By optimising over $N$ for a given value of $\epsilon$,
we will now conclude this analysis with an explicit
non-linear upper bound on this speed.

\begin{theorem}[Robinson Stability] \label{thm:RobinsonStability}

Let $\Omega_\F$ be the enhanced Robinson tiling.
Then $\Omega_\F$ is $f$-stable for $d_B$ on the class of Bernoulli noises $\B$,
with $f(\epsilon)=48 \sqrt[3]{6\epsilon}$.
In particular, $\Omega_\F$ is polynomially stable.

\begin{proof}
To simplify things, we start by bounding
$\left(2^{N+2}+1\right)^2\leq 2^{2N+5}$,
so that we are now trying to minimise
$2\left(4^N\times 2^{9} 3\epsilon+\frac{1}{2^N}\right)$.
If we denote $c(\epsilon)= \sqrt[3]{2^{9} 3\epsilon}= 8\sqrt[3]{3\epsilon}$,
then the upper-bound can be rewritten as $2c
\left( \left(2^Nc\right)^2+\frac{1}{2^Nc}\right)$.

If we treat $x=2^N c$ as a real-valued parameter,
then $x^2+\frac{1}{x}$ is minimal at $x=\sqrt[3]{\frac{1}{2}}$,
equal to $\frac{3}{\sqrt[3]{4}}$.
This gives us a $24\sqrt[3]{6\epsilon}$ bound.
As $N$ must be integer,
we cannot have $N=\log_2\left(\frac{x}{c}\right)$,
but by replacing it with the nearest integer (at distance at most $\frac{1}{2}$),
we obtain the previous bound up to a factor $4^{\frac{1}{2}}=2$,
thus the announced bound.
\end{proof}
\end{theorem}

\bibliographystyle{amsplain}
\bibliography{biblio}
% add below the content of your .bbl file produced by bibtex.

\end{document}